\documentclass[a4paper,11pt,reqno]{amsart}
\usepackage{geometry} 
\usepackage{mathtools}
\usepackage{graphicx}
\usepackage{paralist}
\usepackage[utf8]{inputenc}
\usepackage[T1]{fontenc}
\usepackage{csquotes}
\usepackage[all]{xy}
\usepackage{chngcntr}
\usepackage{amsthm}
\usepackage{amsmath}
\usepackage{amsfonts}
\usepackage{amssymb}
\usepackage{mathrsfs}
\usepackage{MnSymbol}
\usepackage{tikz-cd}
\usepackage{tikz, tikz-cd}
\usepackage{mathdots}
\usepackage{stmaryrd}
\usepackage{cleveref}
\usepackage{fancyhdr}
\usepackage{xcolor}

\newtheorem{thm}{Theorem}[section]
\newtheorem{lem}[thm]{Lemma}
\newtheorem{prop}[thm]{Proposition}
\newtheorem{cor}[thm]{Corollary}
\theoremstyle{remark}

\theoremstyle{claim}

\newtheorem*{thm*}{Theorem}
\newtheorem*{ass}{Bootstrap Assumption}

\theoremstyle{definition}
\newtheorem{definition}[thm]{Definition}

\theoremstyle{Compatibility assumption}
\newtheorem{Compatibility assumption}[thm]{Compatibility condition}

\theoremstyle{example}

\theoremstyle{convention}

\counterwithin{equation}{section}

\theoremstyle{convention}
\newtheorem{rem}[thm]{Remark}

\theoremstyle{convention}
\newtheorem{notation}[thm]{Notation}

\newcommand{\bbR}{\mathbb{R}}

\newcommand{\bbH}{\mathbb{H}}

\calclayout

\begin{document}
\title[Improved decay rate in a stability theorem for hyperbolic metrics]
{Improved decay rate in a stability theorem for hyperbolic metrics}
\author{Frieder J\"ackel}
\thanks
{
The author was supported by the DFG priority program "Geometry at infinity"}
\date{June 13, 2023}

\begin{abstract}
Recently, Ursula Hamenstädt and the author proved a stability result for finite volume hyperbolic metrics in dimension three that does \textit{not} assume any upper volume bounds, but that requires an exponentially fine control of the metric in the thin part of the manifold. We use a bootstrap argument to extend the result allowing for a weaker exponential control of the metric. This is achieved by formulating an abstract axiomatic framework.
\end{abstract}

\maketitle


\section{Introduction}

The construction of Einstein metrics on a closed manifold \(M\) has a long history. Evolving a given metric under Ricci flow, possibly with surgery, has proven to be an extremely fruitful method to construct such metrics. Related to the construction is the stability of Einstein metrics, which adresses the question whether a metric that is almost Einstein (in a suitable sense) is close to an Einstein metric. To prove stability results, a classic approach is to develop a suitable convergence theory of Riemannian manifolds, and use Arzelà-Ascoli type compactness theorems together with arguments by contradiction (see \cite[Theorem 1.1 and Proposition 3.4]{And90}, \cite[Theorem 1.4 and Corollary 1.6]{PW97}, \cite{PetConvSurv}, \cite[Chapter 11]{PetRG}). Alternatively, one can try to use an implicit function theorem applied to the so-called \textit{Einstein operator} to perturb a metric that is almost Einstein (in a suitable sense) to an Einstein metric. This approach, together with a gluing construction, was used by Anderson \cite{And06} and Bamler \cite{Bamler2012} to construct Einstein metrics on manifold obtained from hyperbolic manifolds of finite volume by generalized Dehn filling, or more recently by Fine and Premoselli \cite{FP20} to find negatively curved Einstein metrics on the Gromov-Thurston examples of negatively curved \(4\)-manifolds \cite{GT87}. For a more detailed account, we refer the reader to the introduction in \cite{HJ22}.

Extending earlier work of Tian \cite{Tian}, and using ideas of Bamler \cite{Bamler2012}, Ursula Hamenstädt and the author \cite{HJ22} used the implicit function theorem approach to prove a stability result for negatively curved Einstein of finite volume in dimension three. The goal of this article is to show that this stability result still holds if one of the assumptions is relaxed. Our proof is based on a bootstrap argument. We formulate this bootstrap argument in an abstract axiomatic approach, and hope that this will also be applicable in other contexts. Along the way, we also give an overview of the proof in \cite{HJ22}.

We now state our main result, the generalisation of \cite[Theorem 2]{HJ22}. In its formulation, \(R\) denotes the Riemann curvature endomorphism, and \(M_{\rm small}\) is the \textit{small part} of \(M\), a geometrically defined subset of the thin part of \(M\) (see \Cref{Section - small part}).

 \begin{thm}[Stability of Einstein metrics in dimension $3$]
 \label{Main Theorem} For all \(\alpha \in (0,1)\), \(\Lambda \geq 0\), \(\delta \in (0,2)\), \(b > 1\) and \(\eta > 1\) there
   exist $\varepsilon_0=\varepsilon_0(\alpha,\Lambda,\delta,b,\eta)>0$
   and $C=C(\alpha,\Lambda,\delta,b,\eta) >0$
 with the following property. Let \(M\) be a \(3\)-manifold that admits a complete Riemannian metric \(\bar{g}\) satisfying the following conditions for some \(\varepsilon \leq \varepsilon_0\):
\begin{enumerate}[i)]
\item \({\rm vol}(M,\bar{g}) < \infty\);
\item \(-1-\varepsilon \leq \mathrm{sec}(M,\bar{g}) \leq -1+\varepsilon\);
\item For all \( x \in M_{\rm small}\) it holds
\begin{equation}\label{curvature decay - introduction}
  \max_{\pi \subseteq T_xM}|\mathrm{sec}(\pi)+1|, \,
 | \nabla R|(x), \, | \nabla^2R|(x), \, | \nabla^3R|(x) \leq \varepsilon e^{-\eta d(x,\partial M_{\rm small})};
\end{equation}
\item \(|| {\nabla} {\rm Ric}(\bar{g})||_{C^0(M,\bar{g})} \leq \Lambda\);
\item It holds
\begin{equation}\label{integral condition - introduction}
  e^{bd(x, M_{\rm thick})}\int_Me^{-(2-\delta)r_x(y)}
| {\rm Ric}(\bar{g})+2\bar{g}|_{\bar{g}}^2(y) \, d{\rm vol}_{\bar{g}}(y) \leq \varepsilon^2
\end{equation}
for all \(x \in M\), where $r_x(y)=d(x,y)$.
\end{enumerate}
Then there exists a hyperbolic metric \(g_{\rm hyp}\) on \(M\) so that
\[
	||g_{\rm hyp}-\bar{g}||_{C^{2,\alpha}(M,\bar{g})} \leq C \varepsilon^{1-\alpha}.
\]
Moreover, if additionally \(\bar{g}\) is already hyperbolic outside a region \(U \subseteq M\), and if 
\[
	\int_{U}|{\rm Ric}(\bar{g})+2\bar{g}|_{\bar{g}}^2 \, d{\rm vol}_{\bar{g}} \leq \varepsilon^2,
\]
then for all \(x \in M_{\rm thick}\) it holds
\begin{equation*}\label{exp control away from gluing region}
|g_{\rm hyp}-\bar{g}|_{C^{2,\alpha}}(x) \leq C\varepsilon^{1-\alpha}e^{-(1-\frac{1}{2}\delta)\mathrm{dist}_{\bar{g}}(x,\,U \, \cup \, \partial M_{\rm thick})}.
\end{equation*}
\end{thm}

The improvement of \Cref{Main Theorem}, compared to \cite[Theorem 2]{HJ22} (or more precisely \cite[Theorem 10.1]{HJ22}), lies in relaxing the assumption on the decay rate \(\eta\) in the curvature decay condition (\ref{curvature decay - introduction}) (see the end of \Cref{Section - overview} for a detailed explanation).

We briefly comment on the conditions in \Cref{Main Theorem}. It is important to note that condition \(i)\) only requires the finiteness of the volume, but it does \textit{not} assume any upper bound on the volume.
Consequently, the constants \(\varepsilon_0\) and \(C\) do \textit{not} depend on an upper volume bound. Condition \(ii)\) is natural because in dimension three Einstein metrics have constant sectional curvature. The curvature decay condition (\ref{curvature decay - introduction}) in \(M_{\rm small}\) allows to construct \textit{hyperbolic model metrics} in \(M_{\rm small}\) (see \Cref{Section - small part}) that are used to model the linearisation of the Einstein operator \(\mathcal{L}\) by a linear system of ODEs (see \Cref{Section - averaging operator}). The bound on the covariant derivative of the Ricci tensor in \(iv)\) is necessary for a well-defined notion of Hölder norms (see the proof of \cite[Proposition 2.5]{HJ22} and the references therein). Finally, for the integral condition (\ref{integral condition - introduction}) it is important to note that the exponential weight \(b\) outside the integral is close to \(1\), and the exponential weight \(-(2-\delta)\) inside the integral is close to \(-2\). This makes it possible for the weight inside the integral to absorb the weight outside, so that (\ref{integral condition - introduction}) can be checked without assuming any volume bounds (see \cite[Theorem 11.2]{HJ22} and its proof for an example of this).

An immediate consequence of \Cref{Main Theorem} is the following result, which states that if in dimension three a metric is already hyperbolic except in a region of bounded geometry, then it is close to a hyperbolic metric. Previously known results in this direction had to assume that the entire manifold has bounded geometry (see for example \cite[Proposition 3.4]{And90}, \cite[Corollary 1.6]{PW97}, \cite[Theorem 11.4.16 and Theorem 11.4.17]{PetRG}).

\begin{cor}\label{Main corollary} For any \(\alpha \in (0,1)\), \(\Lambda \geq 0\), \(\iota > 0\) and \(v >0\) there exist \(\varepsilon_0=\varepsilon_0(\alpha,\Lambda, \iota, v)>0\) and \(C=C(\alpha,\Lambda,\iota,v)>0\) with the following property. Let \(M\) be a \(3\)-manifold, and let \(\bar{g}\) be a complete Riemannian metric of finite volume on \(M\) satisfying
\[
	|\sec(M,\bar{g})+1| \leq \varepsilon \quad \text{and} \quad || \nabla {\rm Ric}(\bar{g})||_{C^0(M)} \leq \Lambda.
\]
Assume that there is \(\Omega \subseteq M\) so that
\[
	{\rm inj}(\Omega) \geq \iota, \quad {\rm vol}(\Omega) \leq v \quad \text{and} \quad \sec(M,\bar{g})=-1 \, \text{ outside } \, \Omega.
\]
Then there exists a hyperbolic metric \(g_{\rm hyp}\) on \(M\) satisyfing
\[
	||g_{\rm hyp}-\bar{g}||_{C^{2,\alpha}(M,\bar{g})} \leq C \varepsilon^{1-\alpha}.
\]
\end{cor}

\Cref{Main corollary} already follows from \cite[Theorem 2]{HJ22}. However, since it is not explicitly stated in \cite{HJ22}, we formulate it here for the convenience of the reader.

The article is structured as follows. In \Cref{Section - General Background} we review the relevant results from \cite{HJ22}. More specifically, \Cref{Section - overview} introduces the \textit{Einstein operator}, provides an overview of the proof of \cite[Theorem 2]{HJ22}, and explains that to prove \Cref{Main Theorem} only certain \textit{growth estimates} have to be adjusted. In \Cref{Section - small part} we recall the definition of the small part \(M_{\rm small}\) introduced in \cite{HJ22} and show that the curvature decay assumption (\ref{curvature decay - introduction}) allows to construct hyperbolic model metrics in \(M_{\rm small}\) that will be crucial to model the linearisation \(\mathcal{L}\) of the Einstein operator. How \(\mathcal{L}\) can be modelled by a linear system of ODEs by using an \textit{averaging operator} is explained in \Cref{Section - averaging operator}. In \Cref{Section - growth estimate statement} we state the growth estimates needed to prove \Cref{Main Theorem}. As alluded to, we formulate a general axiomatic approach in \Cref{Section - axiomatic formulation}, and show in \Cref{Section - proof of axiomatic growth estimate} that given the axiomatic set up, one can use a bootstrap argument to obtain growth estimates similar to those needed for the proof of \Cref{Main Theorem}. In \Cref{Section - check compatibility conditions} we verify that the axiomatic conditions are satisfied in the case relevant for \Cref{Main Theorem}. 

\bigskip\noindent
{\bf Acknowledgements:}
I am grateful to Ursula Hamenstädt for helpful comments concerning the presentation of this article.

\section{Background}\label{Section - General Background}

\subsection{Rough overview of the proof}\label{Section - overview}

The goal of this section is to give an overview of the proof of \cite[Theorem 2]{HJ22}. We focus on explaining the main ideas, thus omitting technical details and referring to the relevant sections in \cite{HJ22} for more information. As the complete proof of \cite[Theorem 2]{HJ22} is quite long, we hope that this overview will be of interest independent of the rest of this article.


The overall strategy is to make use of the Inverse Function Theorem to produce a hyperbolic metric. Namely, for every inital metric \(\bar{g}\) there is an operator \(\Phi_{\bar{g}}\), called the \textit{Einstein operator}, that assigns to every metric \(g\) a symmetric \((0,2)\)-tensor \(\Phi_{\bar{g}}(g)\) (see \cite[Section 2.2]{HJ22}). The Einstein operator has two crucial properties.
\begin{enumerate}[(1)]
\item If \(\Phi_{\bar{g}}(g)=0\), then \(g\) is hyperbolic (see \cite[Lemma 2.4]{HJ22}).
\item The linearisation \(\mathcal{L}=(D\Phi_{\bar{g}})_{\bar{g}}\) at the initial metric \(\bar{g}\) is the elliptic partial differential operator acting on symmetric \((0,2)\)-tensors \(h\) by
\[
	\mathcal{L}h=\frac{1}{2}\Delta_L h+ 2h=\frac{1}{2}\Delta h + \frac{1}{2}{\rm Ric}(h)+2h,
\]
where \(\Delta_L\) is the Lichnerowicz Laplacian, \(\Delta\) is the connection Laplacian, and \({\rm Ric}(\cdot)\) is the Weitzenböck curvature operator (see \cite[Section 2.1 and Section 2.2]{HJ22} for more information).
\end{enumerate}
Moreover, one constructs suitable Banach spaces \((S,||\cdot||_{\rm source})\) and \((T,||\cdot||_{\rm target})\) of symmetric \((0,2)\)-tensors, and considers the Einstein operator as a map \(\Phi_{\bar{g}}:U \subseteq S \to T\), where \(U\) is a neighbourhood of \(\bar{g}\) all whose elements are Riemannian metrics. If \(\bar{g}\) is almost hyperbolic (in a suitable sense) and \(||\cdot||_{\rm target}\) is suitably defined, then the image of the inital metric \(||\Phi_{\bar{g}}(\bar{g})||_{\rm target}\) will be close to \(0\). The linearisation \(\mathcal{L}=(D\Phi_{\bar{g}})_{\bar{g}}: S \to T\) is an elliptic operator, and for a suitable choice for the norms \(||\cdot||_{\rm source}\) and \(||\cdot||_{\rm target}\) it will be an invertible operator. Therefore, by the Inverse Function Theorem \(\Phi_{\bar{g}}\) is invertible in a neighbourhood \(V\) of \(\Phi_{\bar{g}}(\bar{g})\), and if \(V\) is large enough so that \(0 \in V\), then there exists a preimage of \(0\). This preimage will be the desired hyperbolic metric by (1). Here we remark that if \({\rm dim}(M)\geq 4\), then \(\Phi_{\bar{g}}(g)=0\) will only imply that \(g\) is an Einstein metric with \({\rm Ric}(g)=-(n-1)g\).

Showing that \(\mathcal{L}:S \to T\) is an invertible operator between the suitably defined Banach spaces \(S\) and \(T\) of symmetric \((0,2)\)-tensors boils down to proving an a priori estimate
\begin{equation}\label{eq: a priori estimate for L}
	||h||_{\rm source} \leq C ||\mathcal{L}h||_{\rm target}
\end{equation}
for suitably constructed norms \(||\cdot||_{\rm source}\) and \(||\cdot||_{\rm target}\). We stress the fact that the constants \(\varepsilon_0\) and \(C\) appearing in \Cref{Main Theorem} have to be \textit{universal constants}, that is, they are only allowed to depend on the stated constants, and in particular can \textit{not} depend on an upper volume bound. Hence also the constant \(C\) in the a priori estimate (\ref{eq: a priori estimate for L}) has to be a universal constant.

A central idea of \cite{HJ22}, which is taken from Tian \cite{Tian}, is to define certain \textit{hybrid norms} that are a mixture of Hölder and (weighted) \(L^2\)-norms in order to guarantee that the a priori estimate (\ref{eq: a priori estimate for L}) can be achieved with a constant \(C\) independent of \({\rm vol}(M)\) (see \cite[Section 4.1]{HJ22} for the precise definition of the hybrid norms). Due to Schauder estimates (\cite[Proposition 2.5]{HJ22})
\[
	||h||_{C^{2,\alpha}(M)} \leq C \Big(||\mathcal{L}h||_{C^{0,\alpha}(M)}+||h||_{C^0(M)} \Big)
\]
bounding the Hölder norm of \(h\) reduces to proving a \(C^0\)-estimate. To further bound the \(C^0\)-norm of \(h\) (independent of \({\rm vol}(M)\)) we invoke the De Giorgi-Nash-Moser estimate (\cite[Lemma 2.8]{HJ22})
\[
	|h|(x_0) \leq C \Big( ||\mathcal{L}h||_{L^2(B(x_0,1))}+||h||_{L^2(B(x_0,1))} \Big),
\]
where \(|h|(x_0)\) is the pointwise \(C^0\)-norm of \(h\) at \(x_0 \in M\), and \(||\cdot||_{L^2(B(x_0,1))}\) is the \(L^2\)-norm in the ball \(B(x_0,1)\).
As one can bound bound \(||h||_{L^2(M)}\) in terms of \(||\mathcal{L}h||_{L^2(M)}\) with some integration by parts argument, this should indicate that it is possible to deduce the a priori estimate (\ref{eq: a priori estimate for L}) with a universal constant \(C\) when defining the norms \(||\cdot||_{\rm source}\) and \(||\cdot||_{\rm target}\) to be hybrid norms that are a mixture of Hölder and \(L^2\)-norms. We refer to the proof of \cite[Proposition 4.3]{HJ22} for details.

The outlined arguments work in regions of the manifold whose injectivity radius is bounded from below by a positive universal constant, but break down in regions where this is not the case. This is because a lower bound on the injectivity radius is required to import elliptic regularity estimates (such as the De Girgio-Nash-Moser estimate) from subsets of \(\bbR^n\) to the manifold. Therefore, other arguments are required to deal with the thin part \(M_{\rm thin}\) of the manifold (the set of points that have no universal lower bound on the injectivity radius - see \Cref{Section - small part}). The idea in \cite{HJ22} was to define the \textit{small part} of the manifold \(M_{\rm small} \subseteq M_{\rm thin}\) (see \Cref{Section - small part}). In \(M_{\rm thin} \setminus M_{\rm small}\) and \(M_{\rm small}\) completely different arguments are employed to obtain the \(C^0\)-estimate.

In \(M_{\rm thin} \setminus M_{\rm small}\) one still uses regularity theory of elliptic operators, but the estimates are applied in the universal cover \(\tilde{M}\) instead of \(M\) itself. Combined with a counting result for the number of preimages (see \cite[Proposition 7.6]{HJ22}), for \(x_0 \in M_{\rm thin} \setminus M_{\rm small}\) this allows to bound the pointwise \(C^0\)-norm \(|h|(x_0)\) in terms of certain weighted \(L^2\)-norms of \(\mathcal{L}h\). These arguments will not be important in the rest of this article, and so we will not elaborate further on them and refer to \cite[Section 7]{HJ22} for more details. As a consequnce, outside the small part we obtain a \(C^0\)-estimate (see (\ref{||f|| 0 lambda}))
\[
	\sup_{M_{\rm thin} \setminus M_{\rm small}}|h| \leq C||\mathcal{L}h||_{0,\lambda,}\]
where \(||\cdot||_{0,\lambda}\) is a norm on the target space of symmetric \((0,2)\)-tensors whose definition is inspired by the pointwise estimates for \(|h|\) in terms of weighted \(L^2\)-norms of \(\mathcal{L}h\). The exact definition of the norms \(||\cdot||_{0,\lambda}\) is not important here, and we refer the reader to \cite[Proposition 7.5 and Definition 9.9]{HJ22} for more details.

To obtain a \(C^0\)-estimate in \(M_{\rm small}\) one exploits the fact that in \(M_{\rm small}\) the linearisation \(\mathcal{L}\) of the Einstein operator can be modelled by an explicit linear system of ODEs. Namely, the main idea, taken from Bamler \cite{Bamler2012}, is to introduce an \textit{averaging operator} that assigns to each tensor \(h\) in \(M_{\rm small}\) another tensor \(\hat{h}\) that only depends on \(r=d(\cdot,\partial M_{\rm small})\) (see \Cref{Section - averaging operator}). This averaging operator commutes (up to small error) with the differential operator \(\mathcal{L}\). So the PDE \(\mathcal{L}h=f\) reduces to the ODE \(\mathcal{L}\hat{h}=\hat{f}\). Using standard ODE estimates we derive certain \textit{growth estimates} for \(\hat{h}\) (and in turn also for \(h\)). The \(C^0\)-estimate in \(M_{\rm small}\) (and also all other a priori estimates in \cite[Section 9.4]{HJ22}) follow from these growth estimates using geometric arguments by contradiction. 

The difference between \Cref{Main Theorem} and \cite[Theorem 2]{HJ22} (or more precisely \cite[Theorem 10.1]{HJ22}) is that in \cite{HJ22} the decay rate \(\eta\) in the curvature decay condition (\ref{curvature decay - introduction}) had to satisfy \(\eta \geq 2+\lambda\) (where \(\lambda \in (0,1)\) is the parameter entering the definition of the norm \(||\cdot||_{0,\lambda}\) - see (\ref{||f|| 0 lambda})), but in this article we only need \(\eta > 1\) (independent of \(\lambda\)). The proof of the growth estimates mentioned in the previous paragraph was the only situation in \cite{HJ22} where the assumption \(\eta \geq 2+\lambda\) was used. Therefore, it suffices to show that the same growth estimates from \cite[Section 9.3]{HJ22} still hold only assuming \(\eta > 1\). The growth estimate inside the small part \(C_{\rm small}\) of a rank 2 cusp \(C\) is contained in \cite[Proposition 9.10]{HJ22}, while \cite[Proposition 9.11]{HJ22} deals with the small part \(T_{\rm small}\) of a Margulis tube \(T\). Because the proofs are very similar, we focus on the case of a rank 2 cusp, that is, we will only present the proof of \cite[Proposition 9.10]{HJ22} under the weaker assumption \(\eta > 1\).

\subsection{The small part, model metrics and trivial Einstein variations}\label{Section - small part}

The goal of this section is to recall the definitions and results from \cite{HJ22} that will be used in the subsequent sections. 

We start with the definition of the \textit{small part} \(M_{\rm small}\). Let \(M\) be a \(3\)-manifold with \(\sec(M) \in [-4,-\frac{1}{4}]\). The \textit{thin part} of \(M\) is
\[
	M_{\rm thin}=\{x \in M \, | \, {\rm inj}_M(x) < \mu\},
\]
where \(\mu\) is a fixed choice for the \textit{Margulis constant} (see \cite[p. 101]{BGS85}). If \(M\) has finite volume, a connected component of the thin part is either rank 2 cusp \(C\) or Margulis tubes \(T\), which is a tubular neighbourhood of a closed geodesic \(\gamma\) whose length is at most \(2\mu\). 

Fix a rank 2 cusp \(C\). Let \(\xi \in \partial_{\infty}\tilde{M}\) be a point corresponding to \(C\), and choose an associated Busemann function \(b_\xi:\tilde{M} \to \bbR\) (see \cite[Section 3]{BGS85} for details about Busemann functions). This induces a Busemann function \(\bar{b}_{\tilde{\xi}}:C \to \bbR\). For \(r \in \bbR\) we denote by \(T(r)\) the level-torus \(T(r):=\{x \in C \, | \, \bar{b}_{\xi}(x)=r\}\). The \textit{small part} of \(C\) is defined as
\[
	C_{\rm small}=\bigcup_{r}\big\{T(r) \, | \, {\rm diam}(T(r))\leq D  \big\},
\]
where \(D\) is a suitably chosen universal constant (see \cite[p. 41]{HJ22} for the definition of \(D\)). The flow of \(-\nabla \bar{b}_{\xi}\) induces a diffeomorphism \(\Psi:\partial C_{\rm small} \times [0,\infty) \xrightarrow{\cong} C_{\rm small}\) so that for \(x=\Psi(y,r) \in C_{\rm small}\) it holds \(d(x,\partial C_{\rm small})=r\). For this reason we often abbreviate \(r(x)=d(x,\partial C_{\rm small})\).

For a Margulis tube \(T\) the definition of \(T_{\rm small}\) is similar, the main difference being that the Busemann function \(\bar{b}_{\xi}\) is replaced by \(r_{\gamma}=d(\cdot,\gamma)\), where \(\gamma\) is the core geodesic of \(T\) (see \cite[Section 7.1]{HJ22} for more details). The \textit{small part} of \(M\) is
\[
	M_{\rm small}=\bigcup_{T}T_{\rm small} \cup \bigcup_{C}C_{\rm small},
\]
where the union is taken over all Margulis tubes \(T\) and rank 2 cusps \(C\). 

As explained at the end of \Cref{Section - overview} we will, for simplicity of presentation, only present the proof of the growth estimate \cite[Proposition 9.10]{HJ22} in the small part of a rank 2 cusp under the weaker assumption \(\eta > 1\), where \(\eta\) is the decay rate in the curvature decay condition (\ref{curvature decay - introduction}). 

Fix a rank 2 cusp \(C\), and assume that the curvature decay condition
\begin{equation}\label{curvatue decay in small part}
	\max_{\pi \subseteq T_xM}|\sec(\pi)+1|, \, |\nabla R|(x), |\nabla^2 R|(x) \leq \varepsilon_0 e^{-\eta d(x, \partial C_{\rm small})} \quad \text{for all }\, x \in C_{\rm small}
\end{equation}
stated in \Cref{Main Theorem} holds. Here \(R\) denotes the Riemann curvature endomorphism. In the decay condition (\ref{curvature decay - introduction}) in \Cref{Main Theorem} we also ask that this decay condition also holds for \(\nabla^3 R\), but for the proof of the growth estimates this is not needed.

It was shown in \cite[Proposition 8.2]{HJ22} that under the curvature decay assumption (\ref{curvatue decay in small part}) there exist certain \textit{hyperbolic model metrics} on \(C_{\rm small}\). The following terminology is used in its formulation. A metric \(g\) on \(T^2 \times I\) (where \(T^2\) is the two torus and \(I\) is an interval) is called a \textit{cusp metric} if it is of the form
\[
	g=e^{-2r}g_{Flat}+dr^2,
\]
where \(g_{Flat}\) is some flat metric on \(T^2\), and \(r\) is the standard coordinate on \(I \subseteq \bbR\). 

\begin{prop}[Proposition 8.2 in \cite{HJ22}]\label{Model metric in cusp} For any \(\eta > 1\) there exists \(\varepsilon_0=\varepsilon_0(\eta)>0\) with the following property. Let \(M\) be a negatively curved \(3\)-manifold of finite volume, and let \(C\) be a rank 2 cusp of \(M\) so that that the curvature decay condition (\ref{curvatue decay in small part}) is satisfied in \(C_{\rm small}\). Then there exists a cusp metric \(g_{cusp}\) on \(C_{\rm small}\) so that for all \(x \in C_{\rm small}\) it holds
\[
	|g-g_{cusp}|_{C^2}(x)=O\big( e^{-\eta r(x)}\big),
\]
where \(r(x)=d(x,\partial C_{\rm small})\).
\end{prop}

Here we use the following convention for the \(O\)-notation. In its formulation \(X\) shall be an arbitrary set, and a constant \(c\) is \textit{universal} if it only depends on the constants stated in \Cref{Main Theorem}.

\begin{notation}\label{big O notation}\normalfont
For functions \(u, \varphi_1,...,\varphi_m: X \to \bbR\) we write \(u=\sum_{k=1}^m O(\varphi_k)\) if there are \textit{universal} constants \(c_k\) such that \(|u(x)| \leq \sum_{k=1}^mc_k \varphi_k(x)\) for all \(x \in X\).
\end{notation}

For the existence of model metrics in the small part of a Margulis tube see \cite[Proposition 8.1]{HJ22}. The main idea for the proof of \cite[Proposition 8.2]{HJ22} is to compare the Jacobi equation in \(M\) with the Jacobi equation in the comparison space \(\bar{M}=\bbH^3\). The curvature decay condition (\ref{curvatue decay in small part}) translates into the fact that the \(C^2\)-norm of the coefficients of the Jacobi equations are exponentially close, and so the desired estimate can be deduced from standard ODE comparison estimates. For the case of a rank 2 cusp we also make use of an effective version of the uniformization theorem (see \cite[Lemma 8.5]{HJ22}) that is needed to detect the correct conformal structure of the level-tori \(T(r)\) as \(r \to \infty\). We refer to \cite[Section 8]{HJ22} for more information. 

Finally, we come to the definition of \textit{trivial Einstein variations}. Let \(g_{cusp}\) be the cusp metric on \(C_{\rm small} \cong \partial C_{\rm small} \times [0,\infty)\) given by \Cref{Model metric in cusp}, and denote by \(g_{Flat}\) the flat metric on \(\partial C_{\rm small}\) so that \(g_{cusp}=e^{-2r}g_{Flat}+dr^2\). In the following definition we take coordinates \((x^1,x^2,r)\) on \(C_{\rm small}\cong \partial C_{\rm small} \times [0,\infty)\), where \((x^1,x^2)\) are flat coordinates on \((\partial C_{\rm small},g_{Flat})\) and \(r\) is the standard coordinate on \([0,\infty)\).

\begin{definition}[Definition 9.4 in \cite{HJ22}]\label{def - trivial Einstein variations} A \((0,2)\)-tensor \(v\) on \(C_{\rm small}\) is a \textit{trivial Einstein variation} if it is of the form
\[
	v=e^{-2r}v_{ij}dx^idx^j,
\]
and if \(\sum_i v_{ii}=0\), that is, if the trace of \(v\) with respect to \(g_{Flat}\) vanishes.
\end{definition}

Here \(dx^i\) and \(dx^j\) are understood to be either \(dx^1\) or \(dx^2\), but not \(dr\). This definition is taken from Bamler \cite[Section 3.2]{Bamler2012}, but we state the definition in coordinates that are better suited for our geometric situation.

\begin{rem}\label{Lv=0 for trivial Einstein}\normalfont For a trivial Einstein variation it holds \(\mathcal{L}_{cusp}v=0\), where \(\mathcal{L}_{cusp}\) is the operator \(\frac{1}{2}\Delta_L h+2h\) with respect to the cusp metric \(g_{cusp}\) given by \Cref{Model metric in cusp}. Indeed, this follows from (\ref{Lh=f in a cusp in coordinates}) below.
\end{rem}

We refer to \cite[Remark 9.13]{HJ22} and the discussion after \cite[Remark 9.5]{HJ22} for explanations why trivial Einstein variations are crucial to consider (also see the comments after \Cref{Check conditions (vii)}).

\subsection{The averaging operator}\label{Section - averaging operator}

The main idea to prove the growth estimate \cite[Proposition 9.10]{HJ22} in the small part \(C_{\rm small}\) of a rank 2 cusp \(C\) is that the linearisation \(\mathcal{L}\) of the Einstein operator can be modelled by a linear system of ODEs. Namely, consider the hyperbolic model metric \(g_{cusp}\) on \(C_{\rm small}\) given by \Cref{Model metric in cusp}. With respect to this model metric we define an \textit{averaging operator} that assigns to each tensor \(h\) another tensor \(\hat{h}\) that only depends on \(r=d(\cdot,\partial C_{\rm small})\). Up to a small error, the averaging operator commutes with the differential operator \(\mathcal{L}\) , so that \(\mathcal{L}h=f\) implies \(\mathcal{L}\hat{h}=\hat{f}\) (up to small error). Since \(\hat{h}\) and \(\hat{f}\) only depend on \(r\), the differential equation \(\mathcal{L}\hat{h}=\hat{f}\) is an ODE. Thus we will obtain growth estimates for \(\hat{h}\) (and in turn for \(h\)) by standard ODE results. It is the goal of this section to introduce the averaging operator and its fundamental properties, and state the necessary ODE results used in the following sections.

Fix a cusp metric \(g\) on \(T^2 \times I\), that is, a metric of the form
\[
	g=e^{-2r}g_{Flat}+dr^2,
\]
where \(g_{Flat}\) is some flat metric on \(T^2\), and \(r\) is the standard coordinate on \(I \subseteq \bbR\). We call a covering \(\varphi:\bbR^2 \times I \to T^2 \times I\) \textit{cusp coordinates} if it is of the form \(\varphi(x^1,x^2,r)=(\psi(x^1,x^2),r)\) for some local isometry \(\psi:\bbR^2 \to (T^2,g_{Flat})\). Moreover, we say that a tensor \(h\) on \(T^2 \times I\) \textit{only depends on} \(r\) if in cusp coordinates the coefficients of \(h\) only depend \(r\). This can be reformulated by saying that \(h\) only depends on \(r\) if \(h\) is invariant under the isometric \(\bbR^2\)-action on \(T^2\times I\) that preserves the level tori \(T(r)=T^2\times \{r\}\).

Now we come to the definition of the \textit{averaging operator}. For a function \(u:T^2 \times I \to \bbR\) we define its average \(\hat{u}:T^2 \times I \to \bbR\) via
\[
	\hat{u}(x):=\frac{1}{{\rm area}(T(r))}\int_{T(r)}u(y) \, d{\rm vol}(y),
\]
where \(r=r(x)\) is the \(I\)-coordinate of \(x\), and \(T(r)=T^2\times \{r\}\) is the level-torus containing \(x\). For a \((0,2)\)-tensor \(h\) its average \(\hat{h}\) is the componentwise average of \(h\), that is, 
\[
	(\hat{h})_{ij}(x):=\widehat{h_{ij}}(x)=\frac{1}{{\rm area}(T(r))}\int_{T(r)}h_{ij}(y) \, d{\rm vol}(y),
\]
where the coefficients are with respect to cusp coordinates. The average for tensors of general type is defined analogously. 

The key properties of the averaging operator are collected in the following lemma. In its formulation \(|h|(x)\) denotes the pointwise \(C^0\)-norm of a tensor \(h\) at the point \(x\).

\begin{lem}[Lemma 9.12 in \cite{HJ22}]\label{Properties of averaging operator} Let \(T^2\times I\) be equipped with a cusp metric. The averaging operation \(\hat{\cdot}\) has the following properties:
\begin{enumerate}[i)]
\item \(\hat{h}\) only depends on \(r\);
\item There is a universal constant \(c>0\) so that 
\[
	|\hat{h}|(x) \leq c \frac{1}{{\rm area}(T(r))}\int_{T(r)}|h|(y) \, d{\rm vol}(y).
\]
In particular,  \(|\hat{h}|(x)\leq  c\max_{T(r(x))}|h|\) for a universal constant \(c\);
\item If \(h\) is of class \(C^1\), then the same holds true for \(\hat{h}\), and \(\widehat{\nabla h}=\nabla \hat{h}\);
\item \(\hat{\cdot}\) commutes with taking the trace,
 that is, \({\rm tr}(\hat{h})=\widehat{{\rm tr}(h)}\);
\item If \(h\) is \(C^1\), then 
\[
	|h-\hat{h}|(x)\leq cDe^{- r(x)}\max_{T(r(x))}|h|_{C^1},
\]
where \(D:={\rm diam}(T^2,g_{Flat})\) and \(c\) is a universal constant.
\end{enumerate}
\end{lem}

It follows from \(iii)\) and \(iv)\) that the averaging operator commutes with the differential operator \(\mathcal{L}=\frac{1}{2}\Delta_L+2{\rm id}\), that is, if \(\mathcal{L}h=f\), then also \(\mathcal{L}\hat{h}=\hat{f}\) (here it is understood that \(\mathcal{L}\) is with respect to the given cusp metric). Property \(v)\) will enable us to deduce bounds for \(h\) from bounds on \(\hat{h}\).

Let \(h\) and \(f\) be \((0,2)\)-tensors that only depend on \(r\). Then the PDE \(\mathcal{L}h=f\) is actually just an ODE. Namely, \(\mathcal{L}h=f\) is equivalent to (see \cite[(9.14)]{HJ22})
\begin{equation}\label{Lh=f in a cusp in coordinates}
	\begin{cases}
	(h_{33})^{\prime \prime}\, \, \, \, \, \, \,-2(h_{33})^{\prime}\,\,\,\,\,-4h_{33}&=-2f_{33}  \\
	(e^{r}h_{i3})^{\prime \prime}\,\,-2(e^{r}h_{i3})^{\prime}-3e^{r}h_{i3}&=-2e^{r}f_{i3} \\
	(e^{2r}h_{ij})^{\prime \prime}-2(e^{2r}h_{ij})^{\prime}&=-2e^{2r}f_{ij}+2\delta_{ij}({\rm tr}(h)-h_{33}),
	\end{cases}
\end{equation}
where \((\cdot)^\prime\) denotes \(\frac{d}{dr}\), and \(\partial_3=\frac{\partial}{\partial r}\). The reason why we state the equations for \(h_{33}\), \(e^{r}h_{i3}\) and \(e^{2r}h_{ij}\), and not \(h_{33},h_{i3},h_{ij}\), is that
\begin{equation}\label{norm in cusp coordinates}
	|h|^2=(h_{33})^2+2\sum_{i=1}^2(e^{r}h_{i3})^2+\sum_{i,j=1}^2(e^{2r}h_{ij})^2.
\end{equation}
Moreover, \({\rm tr}(h)\) satisfies (see \cite[(9.15)]{HJ22})
\begin{equation}\label{PDE for trace as ODE}
	{\rm tr}(h)^{\prime \prime}-2{\rm tr}(h)^\prime-4{\rm tr}(h)=-2{\rm tr}(f).
\end{equation}

From (\ref{Lh=f in a cusp in coordinates}), (\ref{norm in cusp coordinates}), (\ref{PDE for trace as ODE}) one can deduce the following coordinate independent description of trivial Einstein variations (see \Cref{def - trivial Einstein variations}).

\begin{rem}\label{Remark - trivial Einstein variations}\normalfont The trivial Einstein variations \(v\) from \Cref{def - trivial Einstein variations} are exactly the solutions in \(C_{\rm small}\) of \(\mathcal{L}_{cusp}v=0\) so that 
\begin{itemize}
\item \(v\) only depends on \(r\);
\item \(v \in L^2(C_{\rm small})\);
\item \(|v|(r)\) is constant.
\end{itemize}
Here \(\mathcal{L}_{cusp}\) is the operator \(\frac{1}{2}\Delta_L+2{\rm id}\) with respect to the cusp metric \(g_{cusp}\) given by \Cref{Model metric in cusp}, and the norm \(|v|\) is computed with respect to \(g_{cusp}\).
\end{rem}

We will use the two following basic ODE estimates to analyze (\ref{Lh=f in a cusp in coordinates}) and (\ref{PDE for trace as ODE}). In their formulation, \(I\) either denotes  \(\bbR_{\geq 0}\) or an interval of the form \([0,R-1]\) for some \(R \geq 2\) (the case \(I=[0,R-1]\) is needed in the situation of a Margulis tube). Moreover, for a polynomial \(Q=\sum_ma_mX^m\) we write \(Q(\frac{d}{dr})\) for the differential operator \(\sum_ma_m\frac{d^m}{dr^m}\).

\begin{lem}\label{Transfer of exponential rates} Let \(Q \in \bbR[X]\) be a quadratic polynomial with two distinct real roots \(\lambda_1\), \(\lambda_2\). Let \(y:I \to \bbR\) be a solution of the ODE
\[
	Q\left(\frac{d}{dr} \right)(y)=u,
\]
where \(u:I \to \bbR\) is a function satisyfing \(u(r)=\sum_{k=1}^m O\big( \beta_k e^{\mu_k r}\big)\) for some \(\beta_k \in \bbR_{\geq 0}\), and \(\mu_k \in \bbR\setminus \{\lambda_1,\lambda_2\}\). Then
\[
	y(r)=A_1e^{\lambda_1 r}+A_2e^{\lambda_2 r}+\sum_{k=1}^mO\big(\beta_k e^{\mu_k r}\big)
\]
for some constants \(A_1,A_2 \in \bbR\).
\end{lem}

\begin{lem}\label{Transfer of exponential rates - L^1 condition} Let \(Q \in \bbR[X]\) be a quadratic polynomial with two distinct real roots \(\lambda_1\), \(\lambda_2\). Let \(y:I \to \bbR\) be a solution of 
\[
	Q\left(\frac{d}{dr}\right)(y)=u,
\]
where \(u\) satisfies \(|u(r)|\leq e^{a r}\psi(r)\) for some \(a \in \bbR\) and \(\psi \in L^1(\bbR_{\geq 0})\). Then
\[
	y(r)=A_1e^{\lambda_1 r}+A_2e^{\lambda_2 r}+O\big(||\psi||_{L^1(\bbR_{\geq 0})}e^{ar} \big)
\]
for some \(A_1,A_2 \in \bbR\).
\end{lem}

In \Cref{Transfer of exponential rates} and \Cref{Transfer of exponential rates - L^1 condition}, the universal constant absorbed by \(O(...)\) is allowed to depend on \(\lambda_1\), \(\lambda_2\), and \(a\), but \textit{not} on \(R\) (in case \(I=[0,R-1]\)). See \Cref{big O notation} for our convention of the \(O\)-notation.

\begin{proof}[Proof of \Cref{Transfer of exponential rates} and \Cref{Transfer of exponential rates - L^1 condition}] Both of these lemmas follow easily from the explicit integral formulas for solutions of linear ODEs. 
\end{proof}

\section{Growth estimates}\label{Section - Growth estimates}

At the end of \Cref{Section - overview} we explained that in order to prove \Cref{Main Theorem} we only have to prove that the growth estimates from \cite[Section 9.3]{HJ22} still hold under the assumption that the decay rate \(\eta\) in the curvature decay condition (\ref{curvature decay - introduction}) is at least one, and that then all the other arguments from \cite{HJ22} go through without any modifications. In \Cref{Section - growth estimate statement} we give the precise formulation of these growth estimates. As explained in \Cref{Section - overview} and \Cref{Section - averaging operator} the key ingredient to obtain these growth estimates is that, using the averaging operator \(\hat{\cdot}\) from \Cref{Section - averaging operator}, in the small part \(M_{\rm small}\) the linear partial differential operator \(\mathcal{L}\) can be modelled by a linear system of ODEs. In \Cref{Section - axiomatic formulation} we formulate general conditions on an elliptic operator that can be modelled by an ODE that allow to deduce growth estimates analogous to those needed for the proof of \Cref{Main Theorem}. These general growth estimates are proved in \Cref{Section - proof of axiomatic growth estimate}. Finally, in \Cref{Section - check compatibility conditions} we show that the conditions given in \Cref{Section - axiomatic formulation} are satisfied in the specific situation relevant for \Cref{Main Theorem}.

\subsection{Statement of the growth estimate}\label{Section - growth estimate statement}

In this section we state the analog of the growth estimate \cite[Proposition 9.10]{HJ22}. In its formulation the norm \(||\cdot||_{0,\lambda}\) is used. The complete definition (\cite[Definition 9.9]{HJ22}) is not important for our purposes, but we collect two main properties that will be used below. The norm \(||\cdot||_{0,\lambda}\) depends on a parameter \(\lambda \in (0,1)\) and it is defined for all symmetric \((0,2)\)-tensors \(f \in C^{0,\alpha}\big({\rm Sym}^2(T^\ast M) \big)\) of regularity \(C^{0,\alpha}\) (but might be infinite for some \(C^{0,\alpha}\)-tensors). It has the following properties (see \Cref{big O notation} for our convention of the \(O\)-notation):
\begin{enumerate}[(i)]
\item For every rank 2 cusp \(C\) it holds 
\begin{equation}\label{||f|| 0 lambda}
	|f|(x)=O\big(||f||_{0,\lambda}e^{-\lambda r(x)} \big) \quad \text{for all }x \in C_{\rm small}
\end{equation}
where \(|f|(x)\) denotes the \(C^0\)-norm of \(f\) at \(x\), and \(r(x)=d(x,\partial C_{\rm small})\).
\item If \(h \in C^2\big({\rm Sym}^2(T^\ast M) \big) \cap H^1(M)\) is a solution of \(\mathcal{L}h=f\), then it holds
\begin{equation}\label{C^2 bound on d small part}
	\max_{\partial C_{\rm small}}|h|_{C^2}=O(||f||_{0,\lambda}).
\end{equation}
for every rank 2 cusp \(C\).
\end{enumerate}
Point (i) reflects the fact that \(||\cdot||_{0,\lambda}\) is an exponentially weighted Hölder norm (see \cite[(9.2) and (9.8)]{HJ22}). For (ii) we refer to Proposition 7.5, Definition 9.9 and Equation (9.20) in \cite{HJ22}. Moreover, \(C_\lambda^{0,\alpha}\big({\rm Sym}^2(T^*M)\big)\) denotes the Banach space
\[
	C_\lambda^{0,\alpha}\big({\rm Sym}^2(T^*M)\big)=\Big\{ f \in C^{0,\alpha}\big({\rm Sym}^2(T^*M)\big) \, \big| \, ||f||_{0,\lambda} < \infty \Big\}.
\]
If \(M\) is compact \(C_\lambda^{0,\alpha}\big({\rm Sym}^2(T^*M)\big)\) agrees with \(C^{0,\alpha}\big({\rm Sym}^2(T^*M)\big)\) as a set, but the norms are \textit{not} equivalent with universal constants. If \(M\) is non-compact, then \(C_\lambda^{0,\alpha}\big({\rm Sym}^2(T^*M)\big)\) is a proper subset of \(C^{0,\alpha}\big({\rm Sym}^2(T^*M)\big)\).

We now state the generalisation of \cite[Proposition 9.10]{HJ22}. For its formulation we refer to \Cref{Section - small part} for the definition of \(M_{\rm small}\) and trivial Einstein variations, and to \Cref{big O notation} for our convention of the \(O\)-notation. Also recall that \(R\) denotes the Riemann curvature endomorphism.

\begin{prop}[Growth estimate in a cusp]\label{trivial Einstein variation and a priori estimate in a cusp} For all \(\alpha \in (0,1)\), \(\Lambda \geq 0\), \(\lambda \in (0,1)\), and \(\eta > 1\) and there exists \(\varepsilon_0>0\) with the following property. 

Let \(M\) be a finite volume \(3\)-manifold. Assume
\[
	|\sec+1| \leq \varepsilon_0, \quad  ||\nabla {\rm Ric}(g)||_{C^0(M)} \leq \Lambda,
\]
and that for all \(x \in M_{\rm small}\) it holds
\begin{equation}\label{curvature decay - eq in thm}
	\max_{\pi \subseteq T_xM}|\mathrm{sec}(\pi)+1|, \, |\nabla R|(x), \, |\nabla^2R|(x), \, |\nabla^3 R|(x) \leq \varepsilon_0 e^{-\eta d(x,\partial M_{\rm small})}.
\end{equation}
Let \(f \in C_\lambda^{0,\alpha}\big({\rm Sym}^2(T^*M)\big)\), and let \(h \in C^{2}\big({\rm Sym}^2(T^*M)\big)\) with \(||h||_{C^0(M)} < \infty\) be a solution of 
\[
	\mathcal{L}h=f.
\]
Fix a cusp \(C\) of \(M\). Then there exists a unique trivial Einstein variation \(v\) in \(C_{\rm small}\) satisfying
\[
	\sup_{x \in C_{\rm small}}\left(e^{\lambda r(x)}|h-v|(x) \right) < \infty,
\]
and for this trivial Einstein variation it holds
\[
	|v| =O(||f||_{0,\lambda}).
\]
Moreover, if \(||h||_{C^0(M)},||f||_{C^{0,\alpha}(M)} \leq 1\), then for all \(x \in C_{\rm small}\) it holds
\begin{equation}\label{growth estimate in a cusp - eq}
	|h|(x) = O\Big(||f||_{0,\lambda}+e^{-r(x)}\Big)
\end{equation}
and
\begin{equation}\label{exponential growth estimate for h-v in cup}
	e^{\lambda r(x)}|h-v|(x) =  O\Big(||f||_{0,\lambda}+e^{-(1-\lambda)r(x)}\Big),
\end{equation}
where \(r(x)=d(x,\partial C_{\rm small})\).
\end{prop}

The constants absorbed by \(O(...)\) are allowed to depend on the given parameters \(\alpha, \Lambda, \lambda, \eta\).

Recall that the difference to \cite[Proposition 9.10]{HJ22} is that in \cite{HJ22} the decay rate \(\eta\) in the curvature decay condition (\ref{curvature decay - eq in thm}) had to satisfy \(\eta \geq 2+\lambda\), where \(\lambda \in (0,1)\) is the exponential weight in the definition of \(||\cdot||_{0,\lambda}\) (see (\ref{||f|| 0 lambda})), and here we only require \(\eta > 1\) (independent of the parameter \(\lambda \in (0,1)\)).

The curvature decay condition (\ref{curvature decay - eq in thm}) is used to model \(\mathcal{L}\) by an ODE. More precisely, (\ref{curvature decay - eq in thm}) guarantees the existence of a hyperbolic cusp metric modelling the given metric (\Cref{Model metric in cusp}), and with respect to the hyperbolic cusp metric the operator \(\mathcal{L}\) is given by the linear system of ODEs in (\ref{Lh=f in a cusp in coordinates}). For the existence of the hyperbolic model metric we de \textit{not} need the estimate on \(|\nabla^3 R|(x)\) in (\ref{curvature decay - eq in thm}), but we include it to be coherent with the formulation in \Cref{Main Theorem}.

Here we only state the growth estimate in the small part \(C_{\rm small}\) of a rank 2 cusp \(C\). The growth estimate for the small part of a Margulis tube is contained in \cite[Proposition 9.11]{HJ22}. Since the proofs are similar, we only give the proof for the growth estimate in a cusp under the weaker condition \(\eta > 1\).

\subsection{An axiomatic formulation}\label{Section - axiomatic formulation}

We take an axiomatic approach towards the proof of \Cref{trivial Einstein variation and a priori estimate in a cusp}, that is, we formulate general conditions for some elliptic differential operator of second order and some norm that allow to deduce statements analogous to \Cref{trivial Einstein variation and a priori estimate in a cusp}. We believe that this makes the presentation more transparent, and also makes it clear how to apply our arguments in similar situations.

Let \(N\) be a (possibly non-compact) Riemannian manifold with compact boundary \(\partial N\) (possibly with mutiple connected components). We denote by \(r(x)=d(x,\partial N)\) the distance to the boundary. Using a slight abuse of notation, for a number \(r \in \bbR_{\geq 0}\) we denote by \(N(r)\) the level sets
\begin{equation}\label{def of N(r)}
	N(r):= \big\{x \in N \, | \, r(x)=r \big\}.
\end{equation}
For simplicity we assume that \(N(r)\) is a smooth codimension one submanifold for all \(r \in \bbR_{\geq 0}\).

For a parameter \(\sigma \in \bbR\) and tensors \(h\) and \(f\) we denote by
\begin{equation}\label{weighted H^2-norm - axiomatic}
	||h||_{H^2(N,\,\omega_\sigma)}:=\left(\int_M e^{2\sigma r(y)}|h|_{C^2}^2(y) \, d{\rm vol}(y) \right)^{\frac{1}{2}}
\end{equation}
and
\begin{equation}\label{weighted L^2-norm - axiomatic}
	||f||_{L^2(N,\,\omega_\sigma)}:=\left(\int_M e^{2\sigma r(y)}|f|_{C^0}^2(y) \, d{\rm vol}(y) \right)^{\frac{1}{2}}
\end{equation}
exponentially weighted \(H^2\)- and \(L^2\)-norms, where the symbol \(\omega_\sigma\) should indicate that the weight depends on \(\sigma\). Here \(|\cdot|_{C^k}(y):=\sum_{m=0}^k|\nabla^m \cdot|(y)\) denotes the pointwise \(C^k\)-norm at \(y \in N\). For the pointwise \(C^0\)-norm we will often just write \(|\cdot|(y)\) instead of \(|\cdot|_{C^0}(y)\).

Our axiomatic setup consists of five components:
\begin{itemize}
\item A tensor bundle \(\mathcal{T} \to N\) over the given manifold \(N\);
\item An elliptic partial differential operator of second order \(\mathcal{L}\) acting on \(C^2\)-sections of \(\mathcal{T}\);
\item A subspace \(\mathcal{E}\) of the \(C^2\)-sections of \(\mathcal{T}\);
\item A (possibly infinite) norm \(||\cdot||\) on the space of \(C^0\)-sections of \(\mathcal{T}\);
\item A linear operator \(\hat{\cdot}\) acting on \(C^0\)-sections of \(\mathcal{T}\).
\end{itemize}
In the specific situation of \Cref{trivial Einstein variation and a priori estimate in a cusp} these will be
\begin{equation*}
\begin{gathered}
N=C_{\rm small}, \quad \mathcal{T}={\rm Sym}^2(T^\ast M), \quad
 \mathcal{L}=\frac{1}{2}\Delta_L+2{\rm id}, \\ \quad \mathcal{E}=\{\text{trivial Einstein variations in }C_{\rm small}\}, \quad 
||\cdot||=||\cdot||_{0,\lambda}, \\ \text{and }\hat{\cdot} \text{ is the averaging operator from \Cref{Properties of averaging operator}},
\end{gathered}
\end{equation*}
where \(||\cdot||_{0,\lambda}\) is the norm appearing in \Cref{trivial Einstein variation and a priori estimate in a cusp} (see \Cref{Section - growth estimate statement} for more information), and trivial Einstein variations had been defined in \Cref{def - trivial Einstein variations}. 

The space \(\mathcal{E}\) will typically consist of certain fundamental solutions of \(\mathcal{L}\) or the fundamental solutions of a model operator \(\mathcal{L}_{\rm model}\) that is very close to \(\mathcal{L}\) (see \Cref{Remark - trivial Einstein variations}). The norm \(||\cdot||\) can be infinite for some \(C^0\)-tensors (for example when \(||\cdot||=||\cdot||_{C^{0,\alpha}}\)), and we think of \(||\cdot||\) as defining a Banach space of sections whose \(||\cdot||\)-norm is finite. Finally, one should think of \(\hat{\cdot}\) as an averaging operator (see \Cref{Section - averaging operator}).

We need to make certain assumptions on the data \((N,\mathcal{T},\mathcal{L},\mathcal{E},||\cdot||,\hat{\cdot})\). Most importantly, we require that the elliptic operator \(\mathcal{L}\) admits (weighted) Sobolev estimates. Note that because \(N\) has a boundary, to have a well-posed elliptic problem \(\mathcal{L}h=f\), it is neccesary to impose boundary conditions. For this reason, the (weighted) Sobolev estimates involve boundary data \(\max_{\partial N}|h|_{C^2}\). In the situation relevant for \Cref{trivial Einstein variation and a priori estimate in a cusp}, the boundary data \(\max_{\partial N}|h|_{C^2}\) is bounded in terms of \(||f||\) (see (\ref{C^2 bound on d small part})), but we do \textit{not} include this as an assumption in our axiomatic approach. The other key condition is that \(\hat{\cdot}\) has properties that are typical for an averaging operation. Namely, we ask that \(\hat{\cdot}\) satisfies a \textit{Poincaré inequality} (similar to those in \cite[Section 5.8.1]{EvansPDE}) and property \textit{v)} from \Cref{Properties of averaging operator}. The precise assumptions are listed in \Cref{compatibility assumption} below. For its formulation, we again refer to \Cref{big O notation} for our convention of the \(O\)-notation, and recall that \(|\cdot|(x)\) denotes the pointwise \(C^0\)-norm at \(x \in N\). Also, for any \(k \geq 0\), we denote by \(C^k(\mathcal{T})\) the space of \(C^k\)-sections of the tensor bundle \(\mathcal{T} \to N\).

\begin{Compatibility assumption}\label{compatibility assumption}\normalfont Let \((N,\mathcal{T},\mathcal{L},\mathcal{E},||\cdot||,\hat{\cdot})\) be as above. We say that they are \textit{compatible} if the following conditions are satisfied. 
\begin{enumerate}[(i)]
\item There are weighted Sobolev estimates for \(\mathcal{L}\), that is, for all \(h \in C^2(\mathcal{T})\) it holds
\begin{equation}\label{comp ass - weighted a priori estimate}
	||h||_{H^2(N, \, \omega_\sigma)}=O\Big(||\mathcal{L}h||_{L^2(N, \, \omega_\sigma)}+||h||_{L^2(N, \, \omega_\sigma)}+\sup_{\partial N}|h|_{C^2} \Big),
\end{equation}
where \(||\cdot||_{H^2(N, \, \omega_\sigma)}\) and \(||\cdot||_{L^2(N, \, \omega_\sigma)}\) are the weighted integral norms defined in (\ref{weighted H^2-norm - axiomatic}) and (\ref{weighted L^2-norm - axiomatic}).

\item There are (unweighted) \(L^2\)-estimates for \(\mathcal{L}\), that is, for all \(h \in C^2(\mathcal{T})\) it holds
\begin{equation}\label{comp ass - weighted L^2 for h}
	||h||_{L^2(N)}=O\Big(||\mathcal{L}h||_{L^2(N)}+\sup_{\partial N}|h|_{C^2} \Big).
\end{equation}

\item For all \(\sigma\) in a certain range, say \(\sigma \in [0,b]\), and all \(f \in C^0(\mathcal{T})\)  it holds
\begin{equation}\label{comp ass - weighted L^2 for f}
	||f||_{L^2(N, \, \omega_\sigma)}=O(||f||).
\end{equation}

\item For all \(v \in \mathcal{E}\) it holds
\begin{equation}\label{comp ass - regularity for v}
||v||_{C^2(N)}=O(||v||_{C^0(N)}),
\end{equation}
and if \(\sigma \in [0,b]\), then also
\begin{equation}\label{comp ass - weighted L^2 for v}
	||\mathcal{L}v||_{L^2(N, \, \omega_\sigma)}=O(||v||_{C^0(N)}).
\end{equation}
Moreover, the elements of \(\mathcal{E}\) are invariant under \(\hat{\cdot} \,\), that is, 
\begin{equation}\label{comp ass - hat(v)=v}
	\hat{v}=v \quad \text{for all }v \in \mathcal{E}.
\end{equation}

\item The operator \(\hat{\cdot}\) satisfies a level-wise Poincaré inequality, that is, there is a universal constant \(C\) so that for every \(h \in C^1(\mathcal{T})\) and for all \(r \in \bbR_{\geq 0}\) it holds
\begin{equation}\label{comp ass - poincare inequality}
	\int_{N(r)}|h-\hat{h}|^2(y) \, d{\rm vol}(y) \leq C  {\rm diam}(N(r))^2\int_{N(r)}|h|_{C^1}^2(y) \, d{\rm vol}(y),
\end{equation}
where \(N(r)\) is the codimension one submanifold defined in (\ref{def of N(r)}).

\item There is a universal constant \(C\) so that for all \(h \in C^1(\mathcal{T})\) and all \(x \in N\) it holds
\begin{equation}\label{comp ass - ptwise h-hat(h) estimate}
	|h-\hat{h}|(x)\leq C {\rm diam}(N(r(x))) \max_{N(r(x))}|h|_{C^1},
\end{equation}
where \(N(r(x))\) is the level set containing \(x\).

\item For all but finitely many values of \(\sigma\), say \(\sigma \in [0,b]\setminus \{\sigma_1,...,\sigma_k\}\), the following holds. Let \(h\) be a \(C^2\)-tensor with \(h \in L^2(N)\), and assume there exists some \(v \in \mathcal{E}\) so that
\begin{equation}\label{comp ass - bootstrap assumption}
	||v||_{C^0(N)}=O\Big(||\mathcal{L}h||+\max_{\partial N}|h|\Big) \quad \text{and} \quad ||h-v||_{H^2(N, \, \omega_\sigma)}=O\Big(||\mathcal{L}h||+\max_{\partial N}|h|\Big),
\end{equation}
where \(||\cdot||_{H^2(N, \, \omega_\sigma)}\) is the weighted Sobolev norm defined in (\ref{weighted H^2-norm - axiomatic}). Then there exists a (possibly different) \(v^\prime \in \mathcal{E}\) satisfying
\begin{equation}\label{comp ass - v'}
	||v^\prime||_{C^0(N)}=O\Big(||\mathcal{L}h||+\max_{\partial N}|h|\Big) 
\end{equation}
and
\begin{equation}\label{comp ass - hat(h)-v}
	\quad |\hat{h}-v^\prime|_{C^0}(x)=O\Big(\big(||\mathcal{L}h||+\max_{\partial N}|h|\big)e^{\mu(\sigma)r(x)}\Big)  \text{ for all }x \in N,
\end{equation}
where \(\mu(\sigma)\) is some number depending on \(\sigma\) so that
\begin{equation}\label{comp ass - mu(sigma)}
	\int_N e^{2(\sigma^\prime+\mu(\sigma))r(y)} \, d{\rm vol}(y)=O(1)
\end{equation}
for all \(\sigma^\prime\) that are larger than \(\sigma\) by some definite amount, that is, for all \(\sigma^\prime < \sigma + s_0\), where \(s_0 > 0\) is a fixed positive constant (independent of \(\sigma\)).
\end{enumerate}
\end{Compatibility assumption}

The constants absorbed by \(O(...)\) are allowed to depend on \(\sigma\). It will be shown in \Cref{Section - check compatibility conditions} that these conditions hold in the specific situation needed for \Cref{trivial Einstein variation and a priori estimate in a cusp}. 

The key technical condition in \Cref{compatibility assumption} is (vii), stating that if \(h \in C^2(\mathcal{T})\) can be modelled in the \(H^2(N, \,\omega_\sigma)\)-norm by some controlled \(v \in \mathcal{E}\), then there is some \(v^\prime \in \mathcal{E}\) so that the \(C^0\)-norm of \(h-v^\prime\) is exponentially controlled. Roughly speaking, one should be able check this if the elliptic problem \(\mathcal{L}h=f\) can be modelled by an ODE (see the comments after \Cref{Check conditions (vii)} for a detailed discussion). We also recall that if the averaging operator \(\hat{\cdot}\) (almost) commutes with the operator \(\mathcal{L}\), one can use \(\hat{\cdot}\) to model the PDE \(\mathcal{L}h=f\) by the ODE \(\mathcal{L}\hat{h}=\hat{f}\) (see \Cref{Section - averaging operator}).

Let us comment on the remaining conditions in \Cref{compatibility assumption}. We already discussed the (weighted) Sobolev estimates for \(\mathcal{L}\) and the properties of \(\hat{\cdot}\), explaing conditions (i), (ii), (v) and (vi). Condition (iii) will be satisfied when \(||\cdot||\) is an exponentially weighted norm (similar to (\ref{||f|| 0 lambda})) and the weight entering the definition is big enough to absorb the weight \(e^{2\sigma r(y)}\) in the definition (\ref{weighted L^2-norm - axiomatic}) of \(||\cdot||_{L^2(N, \, \omega_\sigma)}\). We mentioned previously that \(\mathcal{E}\) will typically consist of certain fundamental solutions of \(\mathcal{L}\) (or the fundamental solutions of a model operator \(\mathcal{L}_{\rm model}\) that is very close to \(\mathcal{L}\)). In such a situation one should be able to check (\ref{comp ass - regularity for v}) and (\ref{comp ass - weighted L^2 for v}) by applying standard elliptic regularity theory to the operator \(\mathcal{L}\) (or \(\mathcal{L}_{\rm model}\)).

As alluded to at the beginning of this subsection, in any situation where the \Cref{compatibility assumption} is satisfied, one can obtain estimates analogous to those in \Cref{trivial Einstein variation and a priori estimate in a cusp}. This is contained in the following theorem. Recall that \(r(x)=d(x,\partial N)\).

\begin{thm}[General growth estimate]\label{growth estimate - axiomatic} Let \(N\) be a (possibly non-compact) manifold with compact boundary \(\partial N\) (possibly disconnected). Assume that \((N,\mathcal{T},\mathcal{L},\mathcal{E},||\cdot||,\hat{\cdot})\) are compatible in the sense of \Cref{compatibility assumption}.  In addition, assume that the level sets \(N(r)=\{x \in N \, | \, r(x)=r\}\) satisfy
\begin{equation}\label{diameter of level sets - axiomatic}
	{\rm diam}(N(r))=O(e^{-r}).
\end{equation}
Let \(f\) be a \(C^0\)-section of \(\mathcal{T}\) with \(||f|| < \infty\), and let \(h\) be a \(C^2\)-section of \(\mathcal{T}\) with \(h \in L^2(N)\) that solves
\begin{equation}\label{Lh=f - axiomatic}
	\mathcal{L}h=f,
\end{equation}
and whose boundary condition satisfies
\begin{equation}\label{bound on boundary data - axiomatic}
	\max_{\partial N}|h|_{C^2}=O(||f||).
\end{equation}
Then there exists \(v \in \mathcal{E}\) such that
\begin{equation}\label{grwoth estimate - axiomatic eq 1}
	||v||_{C^0(N)}=O(||f||) \quad \text{and} \quad |\hat{h}-v|(x)=O\big(||f||e^{\mu(b)r(x)}\big),
\end{equation}
where \(b\) is the number appearing in conditions (iii), (iv) and (vii) of \Cref{compatibility assumption}, and \(\mu(b)\) is given by (\ref{comp ass - hat(h)-v}). Moreover, if \(||h||_{C^1(N)} < \infty\), then also
\begin{equation}\label{grwoth estimate - axiomatic eq 2}
	|h-v|(x)=O\Big(||f||e^{\mu(b)r(x)} +||h||_{C^1(N)}e^{-r(x)}\Big).
\end{equation}
\end{thm}

We will see in \Cref{Check conditions (i)-(vi)} and \Cref{Check conditions (vii)} that in the situation of \Cref{trivial Einstein variation and a priori estimate in a cusp} the \Cref{compatibility assumption} is satisfied with \(b=2+\lambda-\eta\) and \(\mu(\sigma)=2-\eta-\sigma\), where \(\eta > 1\) is the decay rate in the curvature decay condition (\ref{curvature decay - eq in thm}), and \(\lambda \in (0,1)\) is the exponential weight in the definition of \(||\cdot||_{0,\lambda}\) (see (\ref{||f|| 0 lambda})). Thus \Cref{trivial Einstein variation and a priori estimate in a cusp} will be an immediate consequence of \Cref{growth estimate - axiomatic} once we have checked the conditions (i)-(vii) of \Cref{compatibility assumption} in that specific situation.

The condition (\ref{diameter of level sets - axiomatic}) is taken to fit the specific context of \Cref{trivial Einstein variation and a priori estimate in a cusp} (see (\ref{level tori - diam and area decay})). In a more general setting, where \(N\) is diffeomorphic to \(\partial N \times [0,\infty)\), and the given metric on \(N\) is (very close to) a warped product metric \(\rho(r)^2 g_{\partial N}+dr^2\), (\ref{diameter of level sets - axiomatic})  will be replaced by \({\rm diam}(N(r))=O(\rho(r))\). In that case one will have to adjust the definitions (\ref{weighted H^2-norm - axiomatic}) and (\ref{weighted L^2-norm - axiomatic}), and consider weighted Sobolev spaces whose weights are adapted to \(\rho\). For the clarity of presentation we will not elaborate further on this more general case, and hope that it will be clear how our methods generalize to different settings.

\subsection{Proof of \Cref{growth estimate - axiomatic}}\label{Section - proof of axiomatic growth estimate}

The goal of this section is to explain how the growth estimates of \Cref{growth estimate - axiomatic} follow from \Cref{compatibility assumption}. This will be done by a bootstrap argument. 

Let \((\mathcal{T},\mathcal{L},\mathcal{E},||\cdot||,\hat{\cdot})\) be as in \Cref{compatibility assumption}. Let \(f\) be a \(C^0\)-section of \(\mathcal{T}\) with \(||f|| < \infty\), and let \(h\) be a \(C^2\)-section of \(\mathcal{T}\) with \(h \in L^2(N)\) that solves (\ref{Lh=f - axiomatic}) and whose boundary data satisfies (\ref{bound on boundary data - axiomatic}), that is, 
\begin{equation*}\label{boundary value problem - axiomatic}
	\mathcal{L}h=f \quad \text{and} \quad \max_{\partial N}|h|_{C^2}=O(||f||).
\end{equation*}

The induction assumption of the bootstrap argument, which depends on a parameter \(\sigma \geq 0\), is the following. Note that under the bound on the boundary data (\ref{bound on boundary data - axiomatic}) this is just the assumption (\ref{comp ass - bootstrap assumption}) in condition (vii) of \Cref{compatibility assumption}.

\begin{ass}[Ass$_{\sigma}$]\normalfont 
There exists \(v \in \mathcal{E}\) satisfying
\[
	||v||_{C^0(N)}=O(||f||) \quad \text{and} \quad ||h-v||_{H^2(N, \, \omega_\sigma)}=O\big(||f||\big),
\]
where \(||\cdot||_{H^2(N, \, \omega_\sigma)}\) is the weighted Sobolev norm defined in (\ref{weighted H^2-norm - axiomatic}). Here the universal constant absorbed by \(O(...)\) is allowed to depend on \(\sigma\) (see \Cref{big O notation}).
\end{ass}

For \(\sigma=0\) the bootstrap assumption (Ass$_{\sigma}$) is satisfied with \(v=0\) because of the \(H^2\)- and \(L^2\)-estimates (\ref{comp ass - weighted a priori estimate}), (\ref{comp ass - weighted L^2 for h}), (\ref{comp ass - weighted L^2 for f}), and the boundary condition (\ref{bound on boundary data - axiomatic}).

The induction step of the bootstrap argument is contained in the following proposition. In its formulation, \([0,b]\setminus \{\sigma_1,...,\sigma_k\}\) and \(s_0>0\) are as in condition (vii) of \Cref{compatibility assumption}.

\begin{prop}[Bootstrap Induction]\label{bootstrap - summary}Let \((N,\mathcal{T},\mathcal{L},\mathcal{E},||\cdot||,\hat{\cdot})\) be as in \Cref{growth estimate - axiomatic}. Assume that \({\rm (Ass}_{\sigma}{\rm )}\) holds for some \(\sigma \in [0,b]\setminus \{\sigma_1,...,\sigma_k\}\). Then \({\rm (Ass}_{\sigma^\prime}{\rm )}\) also holds for all \(\sigma^\prime \in [0,b]\) with \(\sigma^\prime < \sigma + \min\{1,s_0\}\).
\end{prop}

We point out that when passing from  \({\rm (Ass}_{\sigma}{\rm )}\) to  \({\rm (Ass}_{\sigma^\prime}{\rm )}\), the \(v \in \mathcal{E}\) given by \({\rm (Ass}_{\sigma}{\rm )}\) might change to some different \(v^\prime \in \mathcal{E}\) in \({\rm (Ass}_{\sigma^\prime}{\rm )}\).

The key condition in \Cref{compatibility assumption} for the bootstrap argument is condition (vii), and we recall that one should be able to check condition (vii) when the operator \(\mathcal{L}\) can be modelled by an ODE.

\begin{proof}The idea of the proof is simple. Let \(v^\prime \in \mathcal{E}\) be given by condition (vii). We split \(h-v^\prime\) as \((h-\hat{h})+(\hat{h}-v^\prime)\), and use the Poincaré inequality (\ref{comp ass - poincare inequality}) to estimate \(||h-\hat{h}||_{L^2(N, \, \omega_{\sigma^\prime})}\) for \(\sigma^\prime \leq \sigma+1\), and (\ref{comp ass - hat(h)-v}) to bound \(||\hat{h}-v^\prime||_{L^2(N, \, \omega_{\sigma^\prime})}\) for \(\sigma^\prime < \sigma+s_0\). Together with (\ref{comp ass - weighted a priori estimate}) this will yield an estimate for \(||h-v^\prime||_{H^2(N, \, \omega_{\sigma^\prime})}\). 

We now give the detailed argument. We start by estimating \(h-\hat{h}\). Let \(v \in \mathcal{E}\) be given by \({\rm (Ass}_{\sigma}{\rm )}\). Note that \(h-\hat{h}=(h-v)-\widehat{h-v}\) because \(\hat{v}=v\) by (\ref{comp ass - hat(v)=v}). Applying the Poincaré inequality (\ref{comp ass - poincare inequality}) to \(h-v\), and recalling that \({\rm diam}(N(r))=O(e^{-r})\) due to (\ref{diameter of level sets - axiomatic}), we deduce
\begin{equation}\label{h-hat(h) over N(r)}
	\int_{N(r)}|h-\hat{h}|^2(y) \, d{\rm vol}(y)\leq C e^{-2r}\int_{N(r)}|h-v|_{C^1}^2(y) \, d{\rm vol(y)}
\end{equation}
for a universal constant \(C\). As \(r(x)=d(x,\partial N)\) is a distance function, the co-area formula states that \(\int_N u \, d{\rm vol}=\int_0^\infty\big(\int_{N(r)}u \, d{\rm vol} \big)\,dr\) for any integrable function \(u\) on \(N\). Thus, multiplying (\ref{h-hat(h) over N(r)}) by \(e^{2(\sigma+1)r}\), and integrating over \(r\) from \(0\) to \(\infty\), we obtain
\[
	\int_N e^{2(\sigma+1)r(y)}|h-\hat{h}|_{C^0}^2(y) \, d{\rm vol}(y) \leq C \int_N e^{2\sigma r(y)}|h-v|_{C^1}^2(y) \, d{\rm vol}(y).
\]
By \({\rm (Ass}_{\sigma}{\rm )}\) the left hand side is bounded by \(O(||f||^2)\). Therefore,
\begin{equation}\label{h-hat(h) - better weighted estimate}
	||h-\hat{h}||_{L^2(N, \, \omega_{\sigma+1})}=O(||f||).
\end{equation}

The assumption (\ref{comp ass - bootstrap assumption}) in condition (vii) is satisfied because of \({\rm (Ass}_{\sigma}{\rm )}\). Let \(v^\prime \in \mathcal{E}\) be so that (\ref{comp ass - v'}) and (\ref{comp ass - hat(h)-v}) hold. So, recalling the boundary condition (\ref{bound on boundary data - axiomatic}), we have \(||v^\prime||_{C^0(N)}=O(||f||)\) and \(|\hat{h}-v^\prime|(x)=O(||f||e^{\mu(\sigma)r(x)})\). Thus,
\[
	\int_N e^{2\sigma^\prime r(y)}|\hat{h}-v^\prime|^2(y) \, d{\rm vol}(y)=O\left(||f||^2\int_N e^{2(\sigma^\prime+\mu(\sigma))r(y)} \, d{\rm vol}(y)\right),
\]
and by the definition (\ref{comp ass - mu(sigma)}) of \(\mu(\sigma)\) it holds \(\int_N e^{2(\sigma^\prime+\mu(\sigma))r(y)} \, d{\rm vol}(y)=O(1)\) when \(\sigma^\prime < \sigma+s_0\). Hence for \(\sigma^\prime < \sigma+s_0\) we have
\begin{equation}\label{hat(h)-v - better weighted estimate}
	||\hat{h}-v^\prime||_{L^2(N, \, \omega_{\sigma^\prime})}=O(||f||).
\end{equation}
Using the triangle inequality we deduce 
\begin{equation}\label{h-v - better weighted estimate}
	||h-v^\prime||_{L^2(N, \, \omega_{\sigma^\prime})}=O(||f||) \quad \text{for }\sigma^\prime < \sigma+ \min\{1,s_0\}
\end{equation}
from (\ref{h-hat(h) - better weighted estimate}) and (\ref{hat(h)-v - better weighted estimate}). 

To bound \(||h-v^\prime||_{H^2(N, \, \omega_{\sigma^\prime})}\), note that for all \(\sigma^\prime \in [0,b]\) it holds
\[
	||\mathcal{L}(h-v^\prime)||_{L^2(N, \, \omega_{\sigma^\prime})} \leq ||f||_{L^2(N, \, \omega_{\sigma^\prime})}+||\mathcal{L}v^\prime||_{L^2(N, \, \omega_{\sigma^\prime})}=O(||f||),
\] 
where we used (\ref{comp ass - weighted L^2 for f}) to bound \(||f||_{L^2(N, \, \omega_{\sigma^\prime})}\), and (\ref{comp ass - weighted L^2 for v}) and \(||v^\prime||_{C^0(N)}=O(||f||)\) to bound \(||\mathcal{L}v^\prime||_{L^2(N, \, \omega_{\sigma^\prime})}\). Moreover, \(\max_{\partial N}|h-v^\prime|_{C^2}=O(||f||)\) because \(\max_{\partial N}|h|_{C^2}=O(||f||)\) by (\ref{bound on boundary data - axiomatic}), and \(\max_{\partial N}|v^\prime|_{C^2}=O(||f||)\) due to (\ref{comp ass - regularity for v}) and \(||v^\prime||_{C^0(N)}=O(||f||)\). Therefore, together with (\ref{h-v - better weighted estimate}) we obtain \(||h-v^\prime||_{H^2(N, \, \omega_{\sigma^\prime})}=O(||f||)\) for all \(\sigma^\prime \in [0,b]\) with \(\sigma^\prime < \sigma+ \min\{1,s_0\}\) by applying (\ref{comp ass - weighted a priori estimate}) to \(h-v^\prime\). Also recall that \(||v^\prime||_{C^0}=O(||f||)\). Thus the estimates in  \({\rm (Ass}_{\sigma^\prime}{\rm )}\) hold with \(v^\prime\). This completes the proof. 
\end{proof}

The proof of \Cref{growth estimate - axiomatic} is now straightforward. 

\begin{proof}[Proof of \Cref{growth estimate - axiomatic}]As mentioned after the definition of \({\rm (Ass}_{\sigma}{\rm )}\), for \(\sigma=0\) the assumption \({\rm (Ass}_{\sigma}{\rm )}\) is satisfied with \(v=0\) due to the \(H^2\)- and \(L^2\)-estimates (\ref{comp ass - weighted a priori estimate}), (\ref{comp ass - weighted L^2 for h}), (\ref{comp ass - weighted L^2 for f}), and the boundary condition (\ref{bound on boundary data - axiomatic}). Applying \Cref{bootstrap - summary} finitely many times, we deduce that \({\rm (Ass}_{\sigma}{\rm )}\) is also satisfied for \(\sigma=b\). Applying condition (vii) for \(\sigma=b\) we obtain some \(v \in \mathcal{E}\) (denoted by \(v^\prime\) in condition (vii)) so that (\ref{comp ass - v'}) and (\ref{comp ass - hat(h)-v}) hold, that is, \(||v||_{C^0(N)}=O(||f||)\) and \(|\hat{h}-v|(x)=O(||f||e^{\mu(b)r(x)})\). This proves (\ref{grwoth estimate - axiomatic eq 1}). If in addition \(||h||_{C^1(N)}<\infty\), then \(|h-\hat{h}|(x)=O(||h||_{C^1(N)}e^{-r(x)})\) due to (\ref{comp ass - ptwise h-hat(h) estimate}) and (\ref{diameter of level sets - axiomatic}). Therefore, the desired estimate (\ref{grwoth estimate - axiomatic eq 2}) for \(|h-v|(x)\) follows from the triangle inequality. This completes the proof.
\end{proof}

\subsection{Checking the compatibility conditions}\label{Section - check compatibility conditions}

The goal of this section is to show that the conditions in \Cref{compatibility assumption} are satisfied in the situation relevant for \Cref{trivial Einstein variation and a priori estimate in a cusp}. For completeness we also include the proofs of \Cref{trivial Einstein variation and a priori estimate in a cusp} and \Cref{Main corollary} at the end.

To deduce \Cref{trivial Einstein variation and a priori estimate in a cusp} from \Cref{growth estimate - axiomatic}, we have to check that the conditions (i)-(vii)  in \Cref{compatibility assumption} are satisfied for
\begin{equation}\label{specific choice of data}
\begin{gathered}
N=C_{\rm small}, \quad \mathcal{T}={\rm Sym}^2(T^\ast M), \quad
 \mathcal{L}=\frac{1}{2}\Delta_L+2{\rm id}, \\ \quad \mathcal{E}=\{\text{trivial Einstein variations in }C_{\rm small}\}, \quad 
||\cdot||=||\cdot||_{0,\lambda}, \\ \text{and }\hat{\cdot} \text{ is the averaging operator from \Cref{Properties of averaging operator}},
\end{gathered}
\end{equation}
Here \(||\cdot||_{0,\lambda}\) is the norm appearing in \Cref{trivial Einstein variation and a priori estimate in a cusp} (see \Cref{Section - growth estimate statement} for more information), trivial Einstein variations had been defined in \Cref{def - trivial Einstein variations}, and the averaging operator \(\hat{\cdot}\) from \Cref{Properties of averaging operator} is defined with respect to the cusp metric \(g_{cusp}\) given by \Cref{Model metric in cusp}. 

We start by checking the conditions (i)-(vi). Recall that \([0,b]\) denotes the allowed range for the weight parameter \(\sigma\), and that \(\lambda \in (0,1)\) is the exponential weight in the definition of \(||\cdot||_{0,\lambda}\) (see (\ref{||f|| 0 lambda})).

\begin{lem}\label{Check conditions (i)-(vi)} Let \((N,\mathcal{T},\mathcal{L},\mathcal{E},||\cdot||,\hat{\cdot})\) be as in (\ref{specific choice of data}). Then the conditions (i)-(vi) in \Cref{compatibility assumption} are satisfied for any \(b < 1+\lambda\).
\end{lem}

Throughout the proof we will use that the level tori 
\[
	T(r)=\{x \in C_{\rm small} \, | \, d(x,\partial C_{\rm small})=r\}
\] 
satisfy 
\begin{equation}\label{level tori - diam and area decay}
	{\rm diam}(T(r))=O(e^{-r}) \quad \text{and} \quad  {\rm area}(T(r))=O(e^{-2r}).
\end{equation}
Indeed, \({\rm diam}(\partial C_{\rm small})\) is bounded by a universal constant due to the definition of \(C_{\rm small}\) (see \Cref{Section - small part}). Thus, the estimates are clear if diameter and area are computed with respect to the cusp metric \(g_{cusp}\) from \Cref{Model metric in cusp} because, by definition, a cusp metric on \(C_{\rm small}=\partial C_{\rm small} \times [0,\infty)\) is of the form \(g_{cusp}=e^{-2r}g_{Flat}+dr^2\) for some flat metric \(g_{Flat}\) on \(\partial C_{\rm small}\). But the given metric \(g\) and \(g_{cusp}\) are uniformly bi-Lipschitz by \Cref{Model metric in cusp}, and thus these estimates also hold for the given metric.

\begin{proof}
\textit{Condition (i).} A detailed account of how to prove weighted \(H^2\)-estimates for \(\mathcal{L}=\frac{1}{2}\Delta_L+2{\rm id}\) using an integration by parts argument is contained in Step 1 of the proof of \cite[Proposition 4.3]{HJ22} (see (4.5) and (4.7) therein). The only differences are that the integral are taken over a complete manifold \(M\) without boundary, and that the weight functions are \(e^{-(2\sqrt{{\rm dim}(M)-2}-\delta)r_x(y)}\), where \(\delta>0\) is a constant and \(r_x(y)=d(x,y)\) for a basepoint \(x \in M\). However, the only property of \(-(2\sqrt{{\rm dim}(M)-2}-\delta)r_x(\cdot)\) that was used is that it is uniformly Lipschitz, and this is also true for \(2\sigma d(\cdot,\partial C_{\rm small})\). Moreover, the boundary terms appearing when invoking integration by parts can all (up to universal constant) be bounded by \(\max_{\partial C_{\rm small}}|h|_{C^2}\) because \({\rm diam}(\partial C_{\rm small})\)  is bounded by a universal constant due to the definition of \(C_{\rm small}\) (see \Cref{Section - small part}). Therefore, the argument presented in \cite{HJ22} goes through with only minor modifications. 

\textit{Condition (ii).} This is contained in the proof of \cite[Lemma 9.16]{HJ22}, or alternatively follows from \cite[Proposition 3.4]{HJ22}. 

\textit{Condition (iii).} Recall that the norm \(||\cdot||_{0,\lambda}\) satisfies \(|f|(x)=O(||f||_{0,\lambda}e^{-r(x)})\) for \(x \in C_{\rm small}\), where \(r(x)=d(x,\partial C_{\rm small})\) (see (\ref{||f|| 0 lambda})). By the co-area formula it holds \(\int_{C_{\rm small}}u \, d{\rm vol}=\int_0^\infty \left(\int_{T(r)}u \, d{\rm vol} \right)\, dr\) for any integrable function \(u\) on \(C_{\rm small}\). Thus, together with (\ref{level tori - diam and area decay}) we obtain
\[
	\int_{C_{\rm small}}e^{2\sigma r(y)}|f|^2(y) \, d{\rm vol}(y)=O\left(||f||_{0,\lambda}^2\int_0^\infty e^{-2r} e^{2\sigma r}e^{-2\lambda r} \, dr\right).	
\]
Note that \(\int_0^\infty e^{-2r} e^{2\sigma r}e^{-2\lambda r} \, dr\) is finite when \(\sigma \in [0,b]\) and \(b < 1+\lambda\). This proves (\ref{comp ass - weighted L^2 for f}).

\textit{Condition (iv).} Let \(v\) be a trivial Einstein variation (see \Cref{def - trivial Einstein variations}). Recall from \Cref{Remark - trivial Einstein variations} that \(\mathcal{L}_{cusp}v=0\), where \(\mathcal{L}_{cusp}\) is the operator \(\frac{1}{2}\Delta_L+2{\rm id}\) with respect to the cusp metric \(g_{cusp}\) given by \Cref{Model metric in cusp}. Thus, it follows from Schauder estimates that  \(||v||_{C^2}=O(||v||_{C^0})\). This proves (\ref{comp ass - regularity for v}). By \Cref{Model metric in cusp} it holds \(|g-g_{cusp}|_{C^2}(x)=O\big(\varepsilon_0 e^{-\eta r(x)}\big)\), where \(g\) is the given metric, and \(\varepsilon_0>0\) and \(\eta > 1\) are the constants in the curvature decay condition (\ref{curvature decay - eq in thm}). Hence \(\big|(\mathcal{L}-\mathcal{L}_{cusp})v\big|(x)=O\big(\varepsilon_0|v|_{C^2}(x)e^{-\eta r(x)}\big)=O\big(\varepsilon_0||v||_{C^0}e^{-\eta r(x)}\big)\). Thus, as \(\eta > 1 > \lambda\), the argument from the proof of condition (iii) also proves (\ref{comp ass - weighted L^2 for v}). Finally, note that a trivial Einstein variation \(v\) only depends on \(r\) (see \Cref{def - trivial Einstein variations}), and thus \(\hat{v}=v\) by the definition of the averaging operator (see \Cref{Section - averaging operator}). So (\ref{comp ass - hat(v)=v}) is also satisfied. 

\textit{Condition (v).} Since the given metric \(g\) and the cusp metric \(g_{cusp}\) from \Cref{Model metric in cusp} are uniformly bi-Lipschitz, we may work with \(g_{cusp}\) instead of \(g\). This has the advantage that the the induced metric on the level tori \(T(r)\) is a flat metric. Inequality (\(\ast \ast\)) on p. 520 of \cite{GromovMetricStructures07} shows that for any flat 2-torus $T^2$ of diameter $1$, we have $\lambda_1(T^2)\geq e^{-2}$. Together with a scaling argument, this implies that if $T^2$ is a flat $2$-torus of 
\({\rm diam}(T^2) \leq 1\), then \(\lambda_1(T^2) \geq \frac{1}{{\rm diam}(T^2)^2}e^{-2}\). Therefore, for any function \(u\) we have the following Poincaré inequality
\begin{equation}\label{Poincare on level torus}
	\int_{T(r)}|u-\hat{u}|^2 \, d{\rm vol} \leq e^2 {\rm diam}(T(r))^2\int_{T(r)}|\nabla u|^2 \, d{\rm vol},
\end{equation}
where \(\hat{u}\) is the average of \(u\) over \(T(r)\) (see \Cref{Section - averaging operator}). Let \(h\) be a \((0,2)\)-tensor. Applying the Poincaré inequality (\ref{Poincare on level torus}) to each coefficient function \(h_{ij}\) of \(h\), and summing over all indices \(i,j\) implies
\[
	\int_{T(r)}|h-\hat{h}|_{C^0}^2 \, d{\rm vol} \leq C{\rm diam}(T(r))^2 \int_{T(r)}|h|_{C^1}^2 \, d{\rm vol}
\]
for a universal constant \(C\). This proves (\ref{comp ass - poincare inequality}).

\textit{Condition (vi).} This is contained in \textit{v)} of \Cref{Properties of averaging operator}.
\end{proof}

The next lemma states that condition (vii) is satisfied in the situation relevant for \Cref{trivial Einstein variation and a priori estimate in a cusp}. For its formulation recall that \([0,b]\setminus \{\sigma_1,...,\sigma_k\}\) denotes the set of all possible values of \(\sigma\) for which condition (vii) can be applied, that \(\mu(\sigma)\) is the exponential rate in (\ref{comp ass - hat(h)-v}), and that \(s_0 >0\) is the definite step size of improvement from \(\sigma^\prime\) to \(\sigma\) in (\ref{comp ass - mu(sigma)}).

\begin{lem}\label{Check conditions (vii)} Let \((N,\mathcal{T},\mathcal{L},\mathcal{E},||\cdot||,\hat{\cdot})\) be as in (\ref{specific choice of data}). Then condition (vii) in \Cref{compatibility assumption} is satisfied with 
\[
[0,b]\setminus \{\sigma_1,...,\sigma_k\}=[0,2+\lambda-\eta]\setminus \{2-\eta\}, \quad \mu(\sigma)=2-\eta-\sigma, \quad s_0=\eta-1,
\]
where \(\eta > 1\) is the decay rate in the curvature decay condition (\ref{curvature decay - eq in thm}), and \(\lambda \in (0,1)\) is the parameter in the norm \(||\cdot||_{0,\lambda}\).
\end{lem}

The core idea of the proof is as follows. Let \(h\) be a solution of \(\mathcal{L}h=f\). Denote by \(\mathcal{L}_{cusp}\) the operator \(\frac{1}{2}\Delta_L+2{\rm id}\) with respect to the cusp metric \(g_{cusp}\) given by \Cref{Model metric in cusp}. Using \Cref{Model metric in cusp} and \Cref{Properties of averaging operator} one can show \(\mathcal{L}_{cusp}\hat{h}=\widehat{f_c}\), where \(f_c:=\mathcal{L}_{cusp}h\) is very close to \(f\). As \(\hat{h}\) and \(\widehat{f_c}\) only depend on \(r\), \(\mathcal{L}_{cusp}\hat{h}=\widehat{f_c}\) is the linear system of ODEs given by (\ref{Lh=f in a cusp in coordinates}). The fundamental solutions of this linear system of ODEs have the growth rates \(e^{(1-\sqrt{5})r},e^{(1+\sqrt{5})r}, e^{-r},e^{3r},1,e^{2r}\). Using \(e^{-2r}|\hat{h}| \in L^1(\bbR_{\geq 0})\) (which follows from \(h \in L^2(C_{\rm small})\) and (\ref{level tori - diam and area decay})) one can eliminate the fundamental solutions with growth rates \(e^{(1+\sqrt{5})r},e^{3r},e^{2r}\). Hence the only relevant fundamental solutions are those with growth rates \(e^{(1-\sqrt{5})r}, e^{-r},1\). Now the trivial Einstein variation \(v^\prime\) in condition (vii) of \Cref{compatibility assumption} can be read off from the linear system of ODEs (\ref{Lh=f in a cusp in coordinates}). Namley, \(v^\prime\) is the part of \(\hat{h}\) that corresponds to the fundamental solution with growth rate \(1\), that is, the fundamental solution with constant norm (see \Cref{Remark - trivial Einstein variations}). So, by \Cref{Transfer of exponential rates} and \Cref{Transfer of exponential rates - L^1 condition}, \(\hat{h}-v^\prime\) is the sum of fundamental solutions with decay rate \(e^{(1-\sqrt{5})r}, e^{-r}\) and an error term with the same decay rate as \(\widehat{f_c}\) (this decay rate has a \(e^{-\lambda r}\) and a \(e^{(2-\eta-\sigma)r}\) summand - see (\ref{growth of hat(f_c)})). 

At this point we make another comment about condition (vii) in \Cref{compatibility assumption} and when it is neccesary to consider a subspace of special tensors \(\mathcal{E}\) as in \Cref{compatibility assumption}. In many situations it is natural and necessary to consider weighted versions of classical norms, such as weighted Hölder norms, when studying a linear elliptic operator \(\mathcal{L}\). However, it may happen that there are special solutions of \(\mathcal{L}v=0\) (or \(\mathcal{L}_{\rm model}v=0\) for a model operator \(\mathcal{L}_{\rm model}\) close to \(\mathcal{L}\)) that are \(L^2\)-integrable but that have constant pointwise \(C^0\)-norm, so that their weighted Hölder norms will be arbitrary large or even infinite (see \Cref{Remark - trivial Einstein variations}). In practice, the set \(\mathcal{E}\) will consist of such special fundamental solutions. These special fundamental solutions make it impossible to prove a priori estimates for the weighted norms. However, it will often be possible to get the weighted estimates if one subtracts these special solutions (such as in the previous paragraph). Indeed, if the operator \(\mathcal{L}\) is modelled by a linear system of ODEs (such as in (\ref{Lh=f in a cusp in coordinates})), one can hope to prove the estimate (\ref{comp ass - hat(h)-v}) for \(\hat{h}-v^\prime\) by arguments similar to those in the previous paragraph.

After these comments we no give the proof of \Cref{Check conditions (vii)}.

\begin{proof}[Proof of \Cref{Check conditions (vii)}]\textbf{Step 1 (Growth rate bound for \(\mathcal{L}_{cusp}\hat{h}\)):}
Let \(h\) be \(C^2\)-tensor with \(h \in L^2(C_{\rm small})\). Abbreviate
\[
	\mathcal{L}h=f,
\]
where \(\mathcal{L}\) is the elliptic operator \(\mathcal{L}h= \frac{1}{2}\Delta_Lh+2h\) with respect to the given metric \(g\) of \(C_{\rm small}\). Let \(g_{cusp}\) be the cusp metric on \(C_{\rm small}\) given by \Cref{Model metric in cusp}, and denote by \(\mathcal{L}_{cusp}\) the elliptic operator \(\mathcal{L}_{cusp}h=\frac{1}{2}\Delta_Lh+2h\) with respect to the metric \(g_{cusp}\). The goal of Step 1 is to show that, under the assumption (\ref{comp ass - bootstrap assumption}) in condition (vii), it holds
\begin{equation}\label{L hat(h)=hat(f)}
	\mathcal{L}_{cusp}\hat{h}=\widehat{f_c}
\end{equation}
for some tensor \(\widehat{f_c}\) that only depends on \(r\) and that satisfies the growth estimate
\begin{equation}\label{growth of hat(f_c)}
	|\widehat{f_c}|(r)=O\Big(\big(||f||_{0,\lambda}+\max_{\partial C_{\rm small}}|h|\big)e^{-\lambda r}+\varepsilon_0 \psi_{\sigma}(r)e^{(2-\eta - \sigma) r} \Big),
\end{equation}
where \(\varepsilon_0>0\) and \(\eta > 1\) are the constants appearing in the curvature decay condition (\ref{curvature decay - eq in thm}), and \(\psi_\sigma:\bbR_{\geq 0} \to \bbR\) is a function with
\begin{equation}\label{L^1-bound on psi-sigma from Ass-sigma}
	||\psi_{\sigma}||_{L^1(\bbR_{\geq 0})}=O\Big(||f||_{0,\lambda}+\max_{\partial C_{\rm small}}|h|\Big).
\end{equation}

To see this, let \(v\) be a trivial Einstein variation so that the assumption (\ref{comp ass - bootstrap assumption}) in condition (vii) is satisfied, that is,
\begin{equation}\label{bootstrap assumption - proof of ODE analysis 1}
	||v||_{C^0}=O\Big(||f||_{0,\lambda}+\max_{\partial C_{\rm small}}|h|\Big),
\end{equation}
and
\begin{equation}\label{bootstrap assumption - proof of ODE analysis 2}
	\left( \int_{C_{\rm small}}e^{2\sigma r(y)}|h-v|_{C^2}^2(y)\, d{\rm vol}(y)\right)^{\frac{1}{2}}=O\Big(||f||_{0,\lambda}+\max_{\partial C_{\rm small}}|h|\Big),
\end{equation}
where \(r(y)=d(y,\partial C_{\rm small})\). It holds 
\[
	\big|(\mathcal{L}-\mathcal{L}_{cusp})(h-v)\big|(x)=O\big(\varepsilon_0 |h-v|_{C^2}(x)e^{-\eta r(x)} \big)
\]
because \(|g-g_{cusp}|_{C^2}(x)=O\big(\varepsilon_0 e^{-\eta r(x)}\big)\) by \Cref{Model metric in cusp}. Recall from the proof of condition (iv) that \(|\mathcal{L}v|(x)=O(\varepsilon_0 ||v||_{C^0}e^{-\eta r(x)})\). So  
\[
	|\mathcal{L}v|(x)=O\Big(\varepsilon_0 \big(||f||_{0,\lambda}+\max_{\partial C_{\rm small}}|h|\big)e^{-\eta r(x)}\Big)
\]
due to (\ref{bootstrap assumption - proof of ODE analysis 1}). We define \(f_c:=\mathcal{L}_{cusp}h\). Recall that \(\hat{\cdot}\) is the averaging operator of \Cref{Properties of averaging operator} with respect to the cusp metric \(g_{cusp}\), and so \(\hat{\cdot}\) commutes with \(\mathcal{L}_{cusp }\) because of \textit{iii)} and \textit{iv)} in \Cref{Properties of averaging operator}. This establishes (\ref{L hat(h)=hat(f)}). Recall \(\mathcal{L}_{cusp}v=0\) from \Cref{Lv=0 for trivial Einstein}, and \(|f|(x)=O(||f||_{0,\lambda}e^{-\lambda r(x)})\) because of (\ref{||f|| 0 lambda}). Combining all these estimates we obtain
\begin{align*}
	|f_c|(x) =& |\mathcal{L}_{cusp}(h-v)|(x) \\
	\leq & |\mathcal{L}h|(x)+|\mathcal{L}v|(x)+\big|(\mathcal{L}-\mathcal{L}_{cusp})(h-v)\big|(x)\\
	=& O\Big(||f||_{0,\lambda}e^{-\lambda r(x)}+\varepsilon_0 \big(||f||_{0,\lambda}+\max_{\partial C_{\rm small}}|h|\big)e^{-\eta r(x)} +\varepsilon_0 |h-v|_{C^2}(x)e^{-\eta r(x)}\Big) \\
	=& O\Big(\big(||f||_{0,\lambda}+\max_{\partial C_{\rm small}}|h|\big)e^{-\lambda r(x)}+\varepsilon_0 |h-v|_{C^2}e^{-\eta r(x)}\Big),
\end{align*}
where in the last line we used \(\eta > 1 > \lambda\). Invoking \(ii)\) of \Cref{Properties of averaging operator} yields
\begin{align*}
	|\widehat{f_c}|(r) =O\left(\Big(||f||_{0,\lambda}+\max_{\partial C_{\rm small}}|h|\Big)e^{-\lambda r}+\varepsilon_0 e^{(2-\eta) r}\int_{T(r)}|h-v|_{C^2}(y) \, d{\rm vol}_2(y) \right),
\end{align*}
where we used that \({\rm area}(T(r))=O(e^{-2r})\) by (\ref{level tori - diam and area decay}). Here we only write \(|\widehat{f_c}|(r)\), and not \(|\widehat{f_c}|(x)\), because the average \(\widehat{f_c}\) only depends on \(r=d(x,\partial C_{\rm small})\) by \(i)\) of \Cref{Properties of averaging operator}. For any \(\sigma \geq 0\) define the function \(\psi_{\sigma}:\bbR_{\geq 0} \to \bbR\) by
\begin{equation}\label{def of psi}
	 \psi_{\sigma}(r):=e^{\sigma r}\int_{T(r)}|h-v|_{C^2}(y) \, d{\rm vol}_2(y),
\end{equation}
where \(v\) is the trivial Einstein variation from (\ref{bootstrap assumption - proof of ODE analysis 2}). With this choice of \(\psi_\sigma\), (\ref{growth of hat(f_c)}) follows from the above estimate for \(|\widehat{f_c}|\).

It remains to check that \(\psi_\sigma\) defined in (\ref{def of psi}) satisfies (\ref{L^1-bound on psi-sigma from Ass-sigma}). To see this, note that
\[
	\int_0^\infty \psi_{\sigma}(r) \, dr=\int_{C_{\rm small}}e^{\sigma r(y)}|h-v|_{C^2}(y) \, d{\rm vol}(y)
\]
due to the co-area formula. By definition of the small part \({\rm diam}(\partial C_{\rm small})\) is bounded by a universal constant, and hence \({\rm vol}(C_{\rm small})\) is also bounded by a universal constant. Therefore, it follows from the Cauchy-Schwarz inequality and (\ref{bootstrap assumption - proof of ODE analysis 2}) that
\[
	\int_0^\infty \psi_{\sigma}(r) \, dr \leq {\rm vol}(C_{\rm small})^{\frac{1}{2}}\left(\int_{C_{\rm small}}e^{2\sigma r(y)}|h-v|_{C^2}^2(y) \, d{\rm vol}(y) \right)^{\frac{1}{2}}=O\Big(||f||_{0,\lambda}+\max_{\partial C_{\rm small}}|h|\Big).
\]
This completes Step 1.

\textbf{Step 2 (ODE Analysis):}
In Step 1 we showed that (see (\ref{L hat(h)=hat(f)}) and (\ref{growth of hat(f_c)}))
\[
	\mathcal{L}_{cusp}\hat{h}=\widehat{f_c} \quad \text{and} \quad  |\widehat{f_c}|(r)=O\Big(\big(||f||_{0,\lambda}+\max_{\partial C_{\rm small}}|h|\big)e^{-\lambda r}+\varepsilon_0 \psi_{\sigma}(r)e^{-(\eta+\sigma-2) r}\Big),
\]
where \(\psi_{\sigma}\) was defined in (\ref{def of psi}). Moreover, as \(\hat{h}\), \(\widehat{f_c}\) only depend on \(r\), the equation \(\mathcal{L}_{cusp}\hat{h}=\widehat{f_c}\) is the linear system of ODEs given
by (\ref{Lh=f in a cusp in coordinates}). Namely, by (\ref{PDE for trace as ODE}), and the first two equations in (\ref{Lh=f in a cusp in coordinates}) we have
\[
	\begin{cases}
	Q_1(\frac{d}{dr})({\rm tr}(\hat{h}))&=-2{\rm tr}(\widehat{f_c}) \\
	Q_1(\frac{d}{dr})(\hat{h}_{33})&=-2(\widehat{f_c})_{33} \\
	Q_2(\frac{d}{dr})(e^r \hat{h}_{i3})&=-2e^r (\widehat{f_c})_{i3} 
	\end{cases}
\]
for some quadratic polynomials \(Q_1\) and \(Q_2\) with roots \(\{1-\sqrt{5},1+\sqrt{5}\}\) and \(\{-1,3\}\). As \(|\widehat{f_c}|\) satisfies the growth estimate (\ref{growth of hat(f_c)}), and since \(-\lambda \nin \{1\pm \sqrt{5},-1,3\}\) we can apply \Cref{Transfer of exponential rates} and
\Cref{Transfer of exponential rates - L^1 condition}.

We know \(||\psi_{\sigma}||_{L^1(\bbR_{\geq 0})}=O(||f||_{0,\lambda})\) due to (\ref{L^1-bound on psi-sigma from Ass-sigma}). Recall from the formulation of condition (vii) and \Cref{Check conditions (vii)} that we only consider \(\sigma \in [0,b]\setminus \{\sigma_1,...,\sigma_k\}=[0,2+\lambda-\eta]\setminus \{2-\eta\}\). In particular, \(-\lambda \leq 2-\eta-\sigma\). Thus we get from \Cref{Transfer of exponential rates} and \Cref{Transfer of exponential rates - L^1 condition} 
\begin{equation}\label{growth in cusp - growth of coefficients 1}
\begin{cases}
	{\rm tr}(\hat{h})(r)&=a_1 e^{(1-\sqrt{5})r}+a_2e^{(1+\sqrt{5})r}+O\big((||f||_{0,\lambda}+\max_{\partial C_{\rm small}}|h|) e^{(2-\eta-\sigma) r} \big); \\
	\hat{h}_{33}(r)&=b_1 e^{(1-\sqrt{5})r}+b_2e^{(1+\sqrt{5})r}+O\big((||f||_{0,\lambda}+\max_{\partial C_{\rm small}}|h|) e^{(2-\eta-\sigma) r} \big); \\
	e^r\hat{h}_{i3}(r)&=c_1^{(i)}e^{-r}\, \, \, \, \, \, \, \, \, \, \, +c_2^{(i)}e^{3r}\,\,\,\,\,\,\,\,\,\,\,+O\big((||f||_{0,\lambda}+\max_{\partial C_{\rm small}}|h|) e^{(2-\eta-\sigma) r} \big);
	\end{cases}
\end{equation}
for some constants \(a_1,a_2,b_1,b_2,c_1^{(i)},c_2^{(i)} \in \bbR\). Note that \(h \in L^2(C_{\rm small}) \subseteq L^1(C_{\rm small})\) since \(C_{\rm small}\) has finite volume, and \({\rm area}(T(r))=O(e^{-2r})\) by (\ref{level tori - diam and area decay}), where \(T(r)\subseteq C_{\rm small}\) is the torus all whose points have distance \(r\) to \(\partial C_{\rm small}\). 
Hence \(e^{-2r}|\hat{h}|(r) \in L^1(\bbR_{\geq 0})\).
In particular, it holds \(e^{-2r} {\rm tr}(\hat{h}) (r), e^{-2r}\hat{h}_{33}(r), e^{-2r}\big(e^{r}\hat{h}_{i3}(r)\big) \in L^1(\bbR_{\geq 0})\)
$(i=1,2)$ because of (\ref{norm in cusp coordinates}), 
and thus 
\[
	a_2=b_2=c_2^{(i)}=0.
\] 
We know \(|\hat{h}|(0)=O\big( \max_{\partial C_{\rm small}}|h|\big)\) by \(ii)\) from \Cref{Properties of averaging operator}. Hence evaluating at \(r=0\) yields
\[
	a_1,b_1,c_1^{(i)}=O\Big(||f||_{0,\lambda} +\max_{\partial C_{\rm small}}|h|\Big).
\]
Note \(1-\sqrt{5} < -1 < -\lambda \leq 2-\eta-\sigma\) because \(\lambda \in (0,1)\) and \(\sigma \in [0,2+\lambda-\eta]\setminus \{2-\eta\}\). Together with the previous estimates we obtain
\begin{equation}\label{growth in cusp - bound on coefficients 1}
	|{\rm tr}(\hat{h})(r)|, |\hat{h}_{33}(r)|, |e^r\hat{h}_{i3}(r)|=O\Big( \big( ||f||_{0,\lambda}+\max_{\partial C_{\rm small}}|h|\big) e^{(2-\eta-\sigma) r}\Big).
\end{equation}
By the last equation in (\ref{Lh=f in a cusp in coordinates}) there is a quadratic polynomial \(Q_3\) with roots \(0\) and \(2\) so that
\begin{align*}
	Q_3\left(\frac{d}{dr}\right)(e^{2r}\hat{h}_{ij})=& 2\delta_{ij}({\rm tr}(\hat{h})-\hat{h}_{33})-2e^{2r}(\widehat{f_c})_{ij} \\
	=& O\Big(\big( ||f||_{0,\lambda}+\max_{\partial C_{\rm small}}|h|\big)e^{(2-\eta-\sigma) r} \Big) \\
	&+O\Big(\big( ||f||_{0,\lambda}+\max_{\partial C_{\rm small}}|h|\big) e^{-\lambda r}+\varepsilon_0\psi_{\sigma}(r)e^{(2-\eta-\sigma) r} \Big),
\end{align*}
where we used the growth rate of \({\rm tr}(\hat{h}), \hat{h}_{33}\) in (\ref{growth in cusp - bound on coefficients 1}), and the one of \(|\hat{f_c}|\) in (\ref{growth of hat(f_c)}). Since we only consider \(\sigma \in [0,2+\lambda-\eta]\setminus \{2-\eta\}\) and \(\eta >1\), it holds \(2-\eta-\sigma \nin \{0,2\}\), and thus we may invoke \Cref{Transfer of exponential rates} and \Cref{Transfer of exponential rates - L^1 condition}. Using (\ref{L^1-bound on psi-sigma from Ass-sigma}) to estimate \(||\psi_{\sigma}||_{L^1(\bbR_{\geq 0})}\), and recalling \(-\lambda \leq 2-\eta-\sigma\), we conclude
\begin{equation}\label{growth in cusp - growth of coefficients 2}
e^{2r}\hat{h}_{ij}=d_1^{(i,j)}+d_2^{(i,j)}e^{2r}+O\Big(\big( ||f||_{0,\lambda}+\max_{\partial C_{\rm small}}|h|\big) e^{(2-\eta-\sigma) r}\Big)
\end{equation}
for some constants \(d_1^{(i,j)},d_2^{(i,j)} \in \bbR\). As before, \(e^{-2r}|\hat{h}|(r) \in L^1(\bbR_{\geq 0})\) implies \(d_2^{(i,j)}=0\). Again, evaluating at \(r=0\) we obtain
\begin{equation}\label{growth in cusp - bound on coefficients 2}
	d_1^{(i,j)}=O\big( ||f||_{0,\lambda}+\max_{\partial C_{\rm small}}|h|\big).
\end{equation}

If \(2-\eta-\sigma > 0\) we set \(v^\prime=0\). Then (\ref{comp ass - v'}) is trivially satisfied. Since \(2-\eta-\sigma > 0\), the desired estimate (\ref{comp ass - hat(h)-v}) on \(|\hat{h}-v^\prime|\) with \(\mu(\sigma)=2-\eta-\sigma\) follows immediately from (\ref{norm in cusp coordinates}), (\ref{growth in cusp - bound on coefficients 1}), (\ref{growth in cusp - growth of coefficients 2}), (\ref{growth in cusp - bound on coefficients 2}), and the fact that \(d_2^{(i,j)}=0\).

If \(2-\eta-\sigma < 0\) we define an Einstein variation \(v^\prime\) in \(C_{\rm small}\) by \(v^\prime_{ij}(r)=d_1^{(i,j)}e^{-2r}\). Note that \({\rm tr}(\hat{h})=\hat{h}_{33}+{\rm tr}(v^\prime)\), and that \({\rm tr}(v^\prime)\) is constant. Moreover, as \(2-\eta-\sigma < 0\) the bound (\ref{growth in cusp - bound on coefficients 1}) yields \(\hat{h}_{33}(r),{\rm tr}(\hat{h})(r) \xrightarrow{r \to \infty}0\). Hence \({\rm tr}(v^\prime)=0\). Therefore, \(v^\prime\) is indeed a trivial Einstein variation (see \Cref{def - trivial Einstein variations}). The desired estimate (\ref{comp ass - v'}) on \(||v^\prime||_{C^0}\) follows from (\ref{growth in cusp - bound on coefficients 2}) and (\ref{norm in cusp coordinates}). Moreover, (\ref{growth in cusp - bound on coefficients 1}), (\ref{growth in cusp - growth of coefficients 2}), the fact that \(d_2^{(i,j)}=0\), and the definition of \(v\) imply
\begin{equation}\label{hat(h)-v exponential estimate}
	|\hat{h}-v^\prime|(r)=O\Big(\big( ||f||_{0,\lambda}+\max_{\partial C_{\rm small}}|h|\big)e^{(2-\eta-\sigma)r}\Big).
\end{equation}
This is the desired estimate (\ref{comp ass - hat(h)-v}) on \(|\hat{h}-v^\prime|\) with \(\mu(\sigma)=2-\eta-\sigma\).

Thus in either case we have completed the proof of (\ref{comp ass - v'}) and (\ref{comp ass - hat(h)-v}). It remains to check (\ref{comp ass - mu(sigma)}) for \(\mu(\sigma)=2-\eta-\sigma\) and the step size \(s_0=\eta-1\), that is, 
\[
	\int_{C_{\rm small}}e^{2(\sigma^\prime+\mu(\sigma))r(y)} \, d{\rm vol}(y)=O(1) \quad \text{for all }\sigma^\prime < \sigma+s_0.
\]
Indeed, it follows from the co-area formula and (\ref{level tori - diam and area decay}) that
\[
	\int_{C_{\rm small}}e^{2(\sigma^\prime+\mu(\sigma))r(y)} \, d{\rm vol}(y)=\int_0^\infty e^{-2r}e^{2(\sigma^\prime+\mu(\sigma))r} \, dr
\]
is finite if \(-1+\sigma^\prime+\mu(\sigma)<0\), that is, if \(\sigma^\prime < 1-\mu(\sigma)=\sigma+ \eta-1\).
\end{proof}

We now show how \Cref{trivial Einstein variation and a priori estimate in a cusp} follows from \Cref{growth estimate - axiomatic}. As explained at the end of \Cref{Section - overview}, once \Cref{trivial Einstein variation and a priori estimate in a cusp} has been established, the remaining arguments from \cite{HJ22} carry over without any modifications, completing the proof of \Cref{Main Theorem}.

\begin{proof}[Proof of \Cref{trivial Einstein variation and a priori estimate in a cusp}] For \((N,\mathcal{T},\mathcal{L},\mathcal{E},||\cdot||,\hat{\cdot})\) as in (\ref{specific choice of data}), all conditions in \Cref{compatibility assumption} are satisfied with \(b=2+\lambda-\eta\) and \(\mu(\sigma)=2-\eta-\sigma\) by \Cref{Check conditions (i)-(vi)} and \Cref{Check conditions (vii)}. The assumption (\ref{diameter of level sets - axiomatic}) in \Cref{growth estimate - axiomatic} also holds because of (\ref{level tori - diam and area decay}). Let \(h\) and \(f\) be as stated in \Cref{trivial Einstein variation and a priori estimate in a cusp}. We know from (\ref{C^2 bound on d small part}) that the boundary condition (\ref{bound on boundary data - axiomatic}) is satisfied. So we can apply \Cref{growth estimate - axiomatic}. It follows from Schauder estimates that \(||h||_{C^1}<\infty\) since by assumption \(||h||_{C^0}<\infty\) and \(||f||_{C^{0,\alpha}} < \infty\). As \(\mu(b)=-\lambda\), (\ref{grwoth estimate - axiomatic eq 1}) and (\ref{grwoth estimate - axiomatic eq 2}) show that there is a trivial Einstein variation \(v\) with
\[
	|v|=O(||f||_{0,\lambda}) \quad \text{and} \quad |h-v|(x)=O\big(||f||_{0,\lambda}e^{-\lambda r(x)}+||h||_{C^1}e^{-r(x)} \big).
\] 
In particular, \(\sup_{x \in C_{\rm small}}e^{\lambda r(x)}|h-v|(x) < \infty\) as \(\lambda \in (0,1)\). Moreover, if \(||h||_{C^0},||f||_{C^{0,\alpha}} \leq 1\), then \(||h||_{C^1}\) is bounded by a universal constant due to Schauder estimates. Thus the above estimate implies (\ref{exponential growth estimate for h-v in cup}). Note that (\ref{growth estimate in a cusp - eq}) follows from (\ref{exponential growth estimate for h-v in cup}) and \(|v|=O(||f||_{0,\lambda})\). Finally, observe that the uniqueness of such a trivial Einstein variation \(v\) is clear because trivial Einstein variations have constant norm (in cusp coordinates). This completes the proof.
\end{proof}

We end with the proof of \Cref{Main corollary}.

\begin{proof}[Proof of \Cref{Main corollary}]Observe that all constants implicitely depend on the choice of a Margulis constant \(\mu\) (see \Cref{Section - small part}). Therefore, if the constants we produce are allowed to depend on a constant \(\iota > 0\), we can without loss of generality assume \(\mu \leq \iota\). 

Fix \(\delta>0\) and \(b > 1\) so that \(2-\delta > b\), and also fix some \(\eta > 1\). Let \((M,\bar{g})\) be as stated in \Cref{Main corollary}, that is, \((M,\bar{g})\) is a complete Riemannian \(3\)-manifold of finite volume with \(|\sec(\bar{g})+1| \leq \varepsilon\) and \(|| \nabla {\rm Ric}(\bar{g})||_{C^0(M)} \leq \Lambda\), and moreover \(\bar{g}\) is already hyperbolic outside a region \(\Omega\) with \({\rm inj}(\Omega) \geq \iota\) and \({\rm vol}(\Omega) \leq v\). The conditions \(i),ii),iv)\) in \Cref{Main Theorem} are satisfied. It holds \(\Omega \subseteq M_{\rm thick}\) because \({\rm inj}(\Omega) \geq \iota\) and \(\mu \leq \iota\). In particular, \(\sec=-1\) in \(M_{\rm thin}\), and so condition \(iii)\) in \Cref{Main Theorem} trivially holds. Therefore, it suffices to check the integral condition (\ref{integral condition - introduction}). From the curvature assumption it follows
\[
	{\rm Ric}(\bar{g})+2\bar{g}=0 \text{ outside }\Omega \quad \text{and} \quad |{\rm Ric}(\bar{g})+2\bar{g}|^2 \leq c\varepsilon^2 \text{ in }\Omega
\]
for a universal constant \(c > 0\). Hence, by using \(\Omega \subseteq M_{\rm thick}\), \(2-\delta > b\) and \({\rm vol}(\Omega) \leq v\), we obtain for all \(x \in M\)
\begin{align*}
	e^{b d(x,M_{\rm thick})}& \int_M e^{-(2-\delta)d(x,y)}|{\rm Ric}(\bar{g})+(n-1)\bar{g}|^2(y) \, d{\rm vol}(y)   &&\\
	& \leq \int_{\Omega}|{\rm Ric}(\bar{g})+(n-1)\bar{g}|^2(y) \, d{\rm vol}(y) &&\\
	& \leq vc\varepsilon^2 && .
\end{align*}
Thus, if \(\varepsilon > 0\) is small enough (depending on \(v\)) the integral condition (\ref{integral condition - introduction}) is satisfied. Therefore, \Cref{Main Theorem} implies the existence of a hyperbolic metric \(g_{\rm hyp}\) close to \(\bar{g}\). This completes the proof.
\end{proof}

\bigskip

\noindent
MATH. INSTITUT DER UNIVERSIT\"AT BONN\\
ENDENICHER ALLEE 60, 53115 BONN, GERMANY\\
e-mail: fjaeckel@math.uni-bonn.de

\end{document}